\magnification\magstep1
 \baselineskip = 18pt
 \def\p{\n{\bf Proof:\ }}
 \def\n{\noindent}
\def\qed{{\hfill{\vrule height7pt width7pt
depth0pt}\par\bigskip}}

 \def\t{\otimes}
 \def\F{F}
 \def\n{\noindent}
 \def\t{\otimes}
 \def\F{F}
\def\ph{\varphi}
%
%
%

\frenchspacing

\newif\ifpagetitre        \pagetitretrue
\newtoks\hautpagetitre    \hautpagetitre={\hfil}
\newtoks\baspagetitre     \baspagetitre={\hfil}
\newtoks\auteurcourant    \auteurcourant={\hfil}
\newtoks\titrecourant     \titrecourant={\hfil}
\newtoks\hautpagegauche   \newtoks\hautpagedroite
\hautpagegauche={\hfil\tensl\the\auteurcourant\hfil}
\hautpagedroite={\hfil\tensl\the\titrecourant\hfil}

\newtoks\baspagegauche
\baspagegauche={\hfil\tenrm\folio\hfil}
\newtoks\baspagedroite
\baspagedroite={\hfil\tenrm\folio\hfil}
\def\w {\widetilde}
\headline={\ifpagetitre\the\hautpagetitre
\else\ifodd\pageno\the\hautpagedroite
\else\the\hautpagegauche\fi\fi}

\footline={\ifpagetitre\the\baspagetitre
\global\pagetitrefalse
\else\ifodd\pageno\the\baspagedroite
\else\the\baspagegauche\fi\fi}

\font\twbf=cmbx12\font\sc=cmcsc10

\font\tenhelv= phvr at 10pt
\font\sevenhelv = phvr at 7pt
\font\fivehelv= phvr at 5pt
\font\ethelv= phvr at 15pt
\font\tengoth=eufm10
\font\sevengoth=eufm7
\font\fivegoth=eufm5
\font\tenbb=msbm10
\font\sevenbb=msbm7
\font\fivebb=msbm5
\newfam\gothfam
\newfam\bbfam
\textfont\gothfam=\tengoth
\scriptfont\gothfam=\sevengoth
\scriptscriptfont\gothfam=\fivegoth

\textfont\bbfam=\tenbb
\scriptfont\bbfam=\sevenbb
\scriptscriptfont\bbfam=\fivebb

\def\rond{{\scriptstyle\circ}}

\def\date{{\the\day}\
\ifcase\month\or Janvier\or \F\'evrier\or Mars\or Avril
\or Mai\or Juin\or Juillet\or Ao\^ut\or Septembre
\or Octobre\or Novembre\or D\'ecembre\fi\ {\the\year}}

    \def\ie{{\it i.e.\/}\ }
 \def\up#1{\raise 1ex\hbox{\sevenrm#1}}
\def\cqfd{\unskip\kern 6pt\penalty 500
\raise -2pt\hbox{\vrule\vbox to 10pt{\hrule width 4pt\vfill
\hrule}\vrule}\par}
\catcode`\@=12
\def\bg{\bigskip\goodbreak}

\def\ref#1&#2&#3&#4&#5\par{\par{\leftskip = 5em{\noindent
\kern-5em\vbox{\hrule height0pt depth0pt width
5em\hbox{\bf[\kern2pt#1\unskip\kern2pt]\enspace}}\kern0pt}
{\sc\ignorespaces#2\unskip},\
{\rm\ignorespaces#3\unskip}\
{\sl\ignorespaces#4\unskip\/}\
{\rm\ignorespaces#5\unskip}\par}}

\def\exo#1{\goodbreak\vskip 12pt plus 20pt minus 2pt%
\line{\noindent\hss\bf%
\uppercase\expandafter{\romannumeral#1}\hss}\nobreak\vskip
12pt }
\def \titre#1\par{\null\vskip
1cm\line{\hss\vbox{\let\\=\cr\twbf\halign%
{\hfil##\hfil\crcr#1\crcr}}\hss}\vskip 1cm}

\let\eps =\varepsilon

\let\ph=\varphi
\def\frac#1#2{{#1\over#2}}

\def\comp{\;{}^{{}_\vert}\!\!\!{\rm C}}

\def\nat{{{\rm I}\!{\rm N}}}

\def\rat{{\rm Q}\kern-.65em {}^{-{}_/}}

\def\N#1{\left\Vert#1\right\Vert}

\let\isspt=\i
\def\i{\ifmmode\infty\else\isspt\fi}

\def\dess#1by#2(#3){\vbox to #2{\hrule width #1
height 0pt depth 0pt\vfill\special{picture #3}}}
\def\Dess#1by#2(#3 scaled#4){{\dimen10=#1\dimen11=#2
\divide\dimen10 by1000\multiply\dimen10 #4
\divide\dimen11 by1000\multiply\dimen11 #4
\vbox to\dimen11{\hrule width\dimen10
height\z@ depth\z@\vfill\special{picture #3 scaled#4}}}}

\def\equipe{\setbox10=\hbox{\dess36mm by 26mm(analyse
scaled 800)}
\setbox12 = \vbox {\hsize =8cm
\centerline{\ethelv EQUIPE D'ANALYSE}
\vskip-4pt
\centerline{\fivehelv URA 754 -- CNRS}
\centerline{\tenhelv Universit\'e Pierre et Marie Curie -
Paris 6}
\vskip8pt
\centerline{\sevenhelv Tour 46 - 0\quad-\quad Bo\^\i te 186
\quad-\quad 4, place
Jussieu\quad-\quad 75252 PARIS
CEDEX 05}
\centerline{\sevenhelv T\'el : (33-1) 44 27 53 49
\quad-\quad T\'el\'ecopie : (33-1) 44 27 25 55\quad-\quad
lana@ccr.jussieu.fr}}
\null\kern-16mm
\line{\kern-10mm\raise-7pt\box10
\hfil\box12\hskip1cm}
\vskip6pt
\line{\kern-9mm\hrulefill}
\vskip3mm
}

 \overfullrule = 0pt

 \centerline{\bf Bounded linear operators}
 \centerline{\bf between $C^*$-algebras}\bigskip
 \centerline{by}\bigskip
 \centerline{Uffe Haagerup}
 \centerline{and}
 \centerline{Gilles Pisier\footnote*{Partially supported by the
 N.S.F.}}
 \vskip.5in
 \n{\bf Plan}

\n Introduction

 \n {  \S 1. Operators between $C^*$-algebras.}

\n { \S 2. Description of $E_n^k$.}

 \n { \S 3. Random series in non-commutative $L_1$-spaces.}

 \n { \S 4. Complements.}
\vskip24pt
 \n {\bf Introduction}

 Let $u\colon \ A\to B$ be a bounded linear operator between two
 $C^*$-algebras $A,B$. The following result was proved in [P1].

 \proclaim Theorem 0.1. There is a numerical constant $K_1$ such that
 for
 all finite sequences $x_1,\ldots, x_n$ in $A$ we have
 $$\leqalignno{&\max\left\{\left\|\left(\sum u(x_i)^*
 u(x_i)\right)^{1/2}\right\|_B, \left\|\left(\sum u(x_i)
 u(x_i)^*\right)^{1/2}\right\|_B\right\}&(0.1)_1\cr
 \le &K_1\|u\| \max\left\{\left\|\left(\sum
 x^*_ix_i\right)^{1/2}\right\|_A,
 \left\|\left(\sum x_ix^*_i\right)^{1/2}\right\|_A\right\}.}$$

A simpler proof  was
given in [H1].
 More recently an other alternate proof  appeared in
[LPP]. In
 this paper we give a sequence of generalizations of this
inequality.

 The above inequality $(0.1)_1$ appears as the case of ``degree one''
 in
 this sequence. The next case of degree 2 seems particularly
 interesting, we
 now formulate it explicitly.

 Let us assume that $A\subset B(H)$ (embedded as a $C^*$-subalgebra)
 for
 some Hilbert space $H$, and similarly that $B\subset B(K)$. Let
 $(a_{ij})$
 be an $n\times n$ matrix of elements of $A$. We define
\vskip2pt
\centerline{$[(a_{ij})]_{(2)} =
\max\left\{\bigg\|(a_{ij})\bigg\|_{M_n(A)},
 \bigg\|(a^*_{ij})\bigg\|_{M_n(A)},
  \left\|\left(\sum_{ij} a^*_{ij}a_{ij}\right)^{1/2}\right\|_A,
 \left\|\left(\sum_{ij}
 a_{ij}a^*_{ij}\right)^{1/2}\right\|_A\right\}$.}
\vskip2pt Then we have

 \proclaim Theorem 0.2. There is a numerical constant $K_2$ such that
 for
 all $n$ and for all $(a_{ij})$ in $M_n(A)$ we have
 $$[(u(a_{ij}))]_{(2)} \le K_2\|u\| [(a_{ij})]_{(2)}.\leqno
(0.1)_2$$

We recall in passing the following identities for
$a_{ij}\in A$ and $a_i\in A$
$$\eqalign{\|(a_{ij})\|_{M_n(A)}  = \sup
\left\{\left|\sum_{ij} \langle
 y_i, a_{ij} x_j\rangle\right|, \quad x_j,y_i\in
H\right.
  \quad \left.\sum \|x_j\|^2 \le 1,\quad \sum\|y_i\|^2 \le
1\right\}, }$$
 and
$$\eqalign{\left\|\left(\sum a^*_ia_i\right)^{1/2}\right\|_A
= \sup
 \left\{\left| \sum \langle y_i, a_ix_0\rangle \right|, \
x_0 \in H,\  y_i
 \in H\right.
  \  \left.\|x_0\| \le 1, \  \sum \|y_i\|^2 \le
1\right\}.}$$
We will denote
 $$[(a_i)]_{(1)} = \max\left\{\left\| \left( \sum
 a^*_ia_i\right)^{1/2}\right\|_A,\quad \left\|\left(
 \sum a_ia^*_i\right)^{1/2}\right\|_A\right\}.\leqno (0.2)$$
 More generally, let us explain the general case of
 "degree $k$" of our
 main
 result. Let $k\ge 1$. Let $n$ be a fixed integer. We will denote $[n]
 =
 \{1,2,\ldots, n\}$. Let $\{a_J\mid J\in [n]^k\}$ be a family of
 elements of
 $A$ indexed by $[n]^k$. Let us denote by $P_k$ the set of all the
 $2^k$
 subsets (including the void set) of $\{1,2,\ldots, k\}$.

\n  For any $\alpha \subset\{1,\ldots, k\}$ we denote by
$\alpha^c$ the
 complement of $\alpha$ and by
 $$\pi_\alpha\colon \ [n]^k \to [n]^\alpha$$
 the canonical projection, i.e.
 $$\forall\ J = (j_1,\ldots, j_k) \in [n]^k
\qquad \pi(J) = (j_i)
 _{i\in\alpha}.$$
 For any $\alpha$ with $\alpha\ne \emptyset$ and $\alpha^c\ne
 \emptyset$ we
 define
 $$\|(a_J)\|_\alpha = \sup \left\{\left| \sum_{J\in [n]^k} \langle
 a_Jx_{\pi_\alpha(J)}, y_{\pi_{\alpha^c}(J)}\rangle\right|\right\}
 \leqno
 (0.3)$$
 where the supremum runs over all families
 $$\{x_\ell \mid \ell \in
 [n]^\alpha\}\quad \hbox{and}\quad \{y_m \mid m\in [n]^{\alpha^c}\}$$
  of elements  of $H$ such that $\sum \|x_\ell\|^2 \le 1$ and $\sum
 \|y_m\|^2 \le 1$. There is an alternate description, we can identify
 $[n]^k$ with $[n]^{\alpha^c} \times [n]^\alpha$ so that $J\in [n]^k$
 is
 identified with $(i,j)$ with $i = \pi_{\alpha ^c}(J)$,
 $j=\pi_\alpha(J)$.
 Then $\|(a_J)\|_\alpha$ is nothing but the norm of the matrix
 $(a_{ij})$
 acting from $\ell_2([n]^\alpha, H)$ into $\ell_2([n]^c, H)$. For
 $\alpha  =
 \emptyset$, this definition extends naturally to
  $$\eqalign{\|(a_J)\|_\emptyset  = \sup
\left\{\left|\sum_{J\in [n]^k}
 \langle a_J x_0, y_J\rangle\right|\right\}
  = \left\|\left(\sum_{J\in [n]^k} a^*_J
a_J\right)^{1/2}\right\|_A}$$
 where the supremum runs over all $x_0 \in H$, $y_J\in H$ such that
 $\|x_0\|
 \le 1$ and $\sum \|y_J\|^2 \le 1$. Similarly, for $\alpha =
 \{1,\ldots,
 k\}$ we set
 $$\|(a_J)\|_\alpha = \left\|\left( \sum a_J
 a^*_J\right)^{1/2}\right\|_A.$$
 We then define
 $$[(a_J)]_{(k)} = \max_{\alpha \in P_k} \{\|(a_J)\|_\alpha\}. \leqno
 (0.4)$$
 We can now state one of our main results.

 \proclaim Theorem 0.k. For each $k\ge 1$, there is a constant $K_k$
 such
 that for any bounded linear operator $u\colon \ A\to B$,
for any $n\ge 1$
 and for any family $\{a_J\mid J\in [n]^k\}$ in $A$ we have
 $$[(u(a_J))]_{(k)} \le K_k\|u\| [(a_J)]_{(k)}.\leqno
(0.1)_k$$
 Moreover, we have $K_k\le 2^{(3k/2)-1}$.

 The proof is essentially in section 1 (it is completed in
section~2).

 We now reformulate this result in a fashion which emphasizes the
 connection
 with the notion of complete boundedness for which we
refer to [Pa].

 Let $A\subset B(H)$ be a $C^*$-algebra embedded as a $C^*$-subalgebra.
 ($H$
 a Hilbert space.) We denote as usual by $M_n$ the set of all $n\times
 n$
 complex matrices (equipped with the norm of the space $B(\ell^n_2)$)
 and by
 $M_n(A)$ the space $M_n\t A$ equipped with its natural $C^*$-norm,
 induced
 by $B(\ell^n_2(H))$. More generally, let $S\subset B({\cal H})$ be any
 closed linear subspace of $B({\cal H})$ (${\cal H}$ is a Hilbert
 space). We
 call $S$ an ``operator space''.

 We denote by $S\t A$ the completion of the linear space $S\t A$
 equipped
 with the norm induced by $B({\cal H}\t_2 H)$ (here ${\cal H}\t _2 H$
 denotes the Hilbert space tensor product
of ${\cal H}$ and $H$).
 We will repeatedly use
the following fact (for a proof see Lemma 1.5 in [DCH]).
Let $K$ be an arbitrary Hilbert space.  Whenever $u:S\to
B(K)$ is completely bounded, the map $I_A \otimes u
:A\otimes S \to A\otimes B(K)$ is bounded and we have
$$\| I_A \otimes u \|_{A\otimes S \to A\otimes B(K)} \le
\|u \|_{cb} . \leqno(0.5)$$
  Clearly
 $S\t A$ is again an operator space embedded into $B({\cal H} \t_2
 H)$.

 For example, we will need to consider a particular embedding of the
 Euclidean space $\ell^n_2 $ into $M_n\oplus  M_n$ as
follows. (We
 equip $M_n\oplus M_n$ with the norm $\|(x,y)\| = \max\{\|x\|,
 \|y\|\}$,
 for which it clearly is an operator space embedded -- say -- into
 $M_{2n}$
 in a block diagonal way.)  We denote by $E_n$ the subspace of $M_n
 \oplus
 M_n$ formed by all the elements of the form
  $$\left(\matrix{x_1\cr \vdots&\quad\bigcirc\cr x_n\cr}\right) \oplus
 \left(\matrix{x_1 &\ldots&x_n\cr&\cr
&\bigcirc\cr}\right)$$
 with $x_1,\ldots, x_n \in {\bf C}$.
 Let $(e_{ij})$ be the usual basis of $M_n$. We denote by
 $$\delta_i = e_{i1} \oplus e_{1i}$$
 the natural basis of $E_n$, (so that the above element can be written
 as
 $\sum x_i\delta_i$.) As a Banach space, $E_n$ is clearly isometric to
 $\ell^n_2$. More precisely, for any $C^*$-algebra $A$ and for any
 $a_1,\ldots, a_n$ in $A$ we have (this known fact is easy to check)
 $$\leqalignno{\left\|\sum \delta_i \t a_i\right\|_{E_n \t A} &= \max
 \left\{\left\|
 \left( \sum a^*_ia_i\right)^{1/2}\right\|, \left\|\left(\sum
 a_ia^*_i\right)^{1/2}\right\|\right\}&(0.6)\cr
 \noalign{\hbox{or equivalently}}
 &= [(a_i)]_{(1)}}$$
 in the preceding notation.

 \n Let us denote by $E^k_n$  the tensor product
 $$E_n \t \cdots \t E_n \qquad (k \hbox{ times}).$$
 Then, Theorem~0.k implies (and is actually equivalent to) the
 following.

 \proclaim Proposition 0.k. For any $u\colon \ A\to B$
 $$\|I_{E^k_n} \t u \| _{E^k_n\t A\to E^k_n \t B} \le
 2^{(3k/2)-1}\|u\|.$$

This proposition is proved in section~1. In section~2 we
extend (0.6) and
compute the norm
 of an element of $E^k_n \t A$  for $k>1$ to deduce
Theorem~0.k
 from Proposition~0.k.

 In section~3, we develop the viewpoint of [LPP] which dualizes
 inequalities
 such as $(0.1)_1$ or $(0.1)_k$ to compute (an equivalent of) the norm
 of
 certain random series with coefficients  in a non-commutative
 $L_1$-space.
 Let $(\varepsilon_j)_{j\in {\bf N}}$ be an i.i.d. sequence of random
 variables each distributed uniformly over the unimodular complex
 numbers.
 (Such variables are sometimes called Steinhaus variables.) Let $A_*$
 be a
 non-commutative $L_1$-space. Roughly, while [LPP] treats the case of
 $A_*$-valued random variables which depend {\it linearly\/} on the
 sequence
 $(\varepsilon_j)$, we can treat variables which depend bilinearly or
 multilinearly in the variables $(\varepsilon_j)$. For a precise
 statement
 see Theorem~3.6 below.

 It might be useful for some readers to emphasize that the variables
 $(\varepsilon_j)$ can be replaced by independent choices of signs or
 more
 importantly by i.i.d. {\it Gaussian\/} variables. All our results
 remain
 true in this setting, but with different numerical constant, this
 follows
 from the fact (due to N.~Tomczak-Jaegermann) that $A_*$ is of cotype
 2,
 see e.g. [P3] p.~36 for more details. We also would like to draw the
 reader's attention to Kwapie\'n's paper [K] which contains
 ``decoupling
 inequalities'' quite relevant to the situation considered in
 Theorem~3.6
 below.
 Using [K] one can deduce from (3.1) below some ``non-decoupled''
 inequalities. For instance, we can find an equivalent of integrals of
 the
 form $\int\left\|\sum\limits_{1\le i <j\le n} \varepsilon_i
 \varepsilon_j
 x_{ij}\right\|_{A_*} dP$ where $x_{ij} \in A_*$ and
 $(\varepsilon_j)_{j\ge
 1}$ is an i.i.d. sequence of symmetric $\pm 1$ valued random variables
 on a
 probability space $(\Omega,  P)$, and similarly in the multilinear
 case. We
 will not spell out the details.

The results of the first three sections of this paper
rely heavily on the following factorization result proved
in section 1:  The identity map $I_{E_n}$ on the operator
space $E_n$ has a completely bounded factorization
through the von Neumann algebra $VN(F_n)$ associated
with the left regular representation of  the free group
with $n$ generators, i.e. there are $w_n:E_n\to VN(F_n)$
and $v_n:VN(F_n)\to E_n$ such that
$$I_{E_n}=v_n w_n \quad {\rm and }\quad \|v_n\|_{cb}
\|w_n\|_{cb}\le 2.$$
In section 4, we show that for any sequence of
factorizations $I_{E_n}=v_n w_n$ ($n=1,2,...$) of the
identity maps   $I_{E_n}$ through {\it injective} von
Neumann algebras  we have
$$\lim_{n\to \infty}  \|v_n\|_{cb}
\|w_n\|_{cb}=+\infty.$$
Combining these two facts about the factorization
of  $I_{E_n}$
with Voiculescu's recent result ([V1]) that the algebra of
all $n\times n$ matrices over $ VN(F_\infty)$ is isomorphic
(as a  von Neumann algebra)
to $ VN(F_\infty)$, we show at the end of section 4 that
the von Neumann algebra $VN(F_n)$ is not a complemented
subspace of $B(H)$ for any $n\ge 2$. (For very recent
results on similar questions, see [P4,CS].) We also
include several general remarks about the relation between
the existence of a {\it-completely bounded} linear
projection from $B(H)$ onto a subspace $S$ and that of a
{\it  bounded}
 linear projection from $B(\ell_2)\otimes B(H)$
onto $B(\ell_2)\otimes S$. For instance, if $S$ is weak-$*$
closed and if $B(\ell_2)\overline{\t} S$ denotes the
weak-$*$ closure of $B(\ell_2){\t} S$ in $B(\ell_2\t H)$,
  we show that   there is a bounded linear
projection from $B(\ell_2\t H)$ onto
$B(\ell_2)\overline{\t} S$ if and only if there is a
completely bounded one from $B(H)$ onto S.

Finally, we compare the space $E_n$ with
the linear span
$S_n$
of a  free system of  random variables
$\{x_1,...,x_n\}$ in a
$C^*$-probability space $(A,\ph)$ in the sense of
Voiculescu [V1,2]. In particular, in the case of
a  semicircular (or circular) system in Voiculescu's
sense,  we show that there is an isomorphism $u$ from
 $E_n$ onto
 the
operator space $S_n$  such that
$$\|u\|_{cb}\|u^{-1}\|_{cb}\le 2.$$

 \magnification\magstep1
 \baselineskip = 18pt
 \def\n{\noindent}
 \def\t{\otimes}
 \def\F{F}
 \def\qed{{\hfill{\vrule height7pt width7pt
depth0pt}\par\bigskip}}
 \overfullrule = 0pt
\vfill\eject

 \n {\bf \S 1. Operators between $C^*$-algebras.}

 We will use repeatedly the following fact which has been known to the
 first
 author for some time. The main point ((1.2) below) is a refinement of
 one
 of the inequalities of [H2]. (We remind the reader that
we denote simply  by ${C^*_\lambda(\F_n)\t A}$ the minimal
or spatial tensor product which is often denoted by
$C^*_\lambda(\F_n)\t_{\rm
 min} A$.)

 \proclaim Proposition 1.1. Let $\F_n$ denote the free group on $n$
 generators $g_1,\ldots, g_n$, and let $C^*_\lambda(\F_n)$ be the
 reduced
 $C^*$-algebra of $\F_n$, i.e. the $C^*$-algebra generated by the left
 regular representation
 $\lambda\colon \ \F_n \to B(\ell^2(\F_n))$. Then
 \item{(1)} For any $C^*$-algebra $A$ and for any set $(a_g)_{g\in S}$
 of
 elements of $A$ indexed by a finite subset $S$ of
$\F_n$:
 $$\left\|\sum_{g\in S} \lambda(g) \t a_g\right\|_
{C^*_\lambda(\F_n)\t
   A} \ge \max\left\{\left\| \sum_{g\in S} a^*_g
a_g\right\|^{1/2},
 \left\|
\sum_{g\in S} a_ga^*_g\right\|^{1/2}\right\}.\leqno
(1.1)$$
 \item{(2)}  For any $C^*$-algebra $A$ and for any set $(a_g)_{g\in G}$
 of
 elements of $A$ indexed by a subset $S$ of $\{g_1,\ldots, g_n,
 g^{-1}_1,\ldots, g^{-1}_n\}$:
 $$\left\|\sum_{g\in S} \lambda(g) \t a_g\right\|_{C^*_\lambda
(\F_n)\t A} \le 2\max\left\{\left\| \sum_{g\in S}
a^*_ga_g\right\|^{1/2},
 \left\| \sum_{g\in S} a_ga^*_g\right\|^{1/2}\right\}.\leqno (1.2)$$

 \n {\bf Proof.} (1)\ let $(\delta_g)_{g\in G}$ be the standard basis
 of
 $\ell^2(\F_n)$. We may assume that $A \subset B(K)$ for some Hilbert
 space
 $K$. Since the min-tensor product coincide with the spatial tensor
 product,
 we have for all unit vectors $\xi \in K$:
 $$\eqalign{\left\|\sum \lambda(g) \t a_g\right\|
_{C^*_\lambda(\F_n) \t A} &\ge \left\|\sum_{g\in S}
(\lambda(g) \t a_g) (\delta_e \t
 \xi)\right\|\cr
 &= \left\|\sum_{g\in S} \delta _g\t a_g\xi\right\|\cr
 &= \left(\sum_{g\in G} \|a_g\xi\|^2\right)^{1/2}\cr
 &= \left(\left(\sum_{g\in G} a^*_ga_g\right)\xi, \xi\right)^{1/2}}$$
 Taking supremum over all unit vectors $\xi \in K$ we get
 $$\left\| \sum_{g\in S} \lambda(g) \t a_g\right\| _{C^*_\lambda (\F_n)
 \t A} \ge \left\|\sum_{g\in S} a^*_g a_g\right\|^{1/2}.$$
 The same argument applied to the norm of $(\lambda(g) \t a_g)^* =
 \lambda(g^{-1}) \t a^*_g$ gives
 $$\left\|\sum_{g\in G} \lambda(g) \t
a_g\right\|_{C^*_\lambda(\F_n)\t A} \ge
\left\|\sum_{g\in S} a_ga^*_g\right\|^{1/2}.$$
 This proves (1). Note that the statement (1) actually holds in
 $C^*_\lambda(\Gamma)\t A$ for any discrete group $\Gamma$.

 \n (2)\ Consider first the case $S =\{g_1,\ldots, g_n,
 g^{-1}_1,\ldots,
 g^{-1}_n\}$. We can write $\F_n$ as a disjoint union:
 $$\F_n = \{e\} \cup \left\{ \bigcup^n_{i=1} \Gamma^+_i\right\} \cup
 \left\{
 \bigcup^n_{i=1} \Gamma^-_i\right\}$$
 where
 $$\eqalign{\Gamma^+_i = &\hbox{ set of
 reduced words starting with
a  positive power of } g_i,\cr
 \Gamma^-_i = &\hbox{ set of reduced words starting with
a negative power of } g_i.}$$
 Let $e_0, e^+_i$ and $e^-_i$ denote the orthogonal projection of
 $\ell^2(\F_n)$ onto the subspaces ${\bf C} \delta_e$,
 $\ell^2(\Gamma^+_i)$
 and $\ell^2(\Gamma^-_i)$ respectively. Then these projections are
 pairwise
 orthogonal and
 $$e_0  + \sum^n_{i=1} e^+_i+ \sum ^n_{i=1} e^-_i = I_{\ell^2(\F_n)}.$$
 For any $g\in G$ and for any generator $g_i$, the length of the
 reduced
 word for $g_ig$ is either
 $$|g_ig| = |g|+1\quad \hbox{or}\quad |g_ig| = |g|-1.$$
 The first case exactly occurs when $g_ig$ starts with an
element of
 $\Gamma^+_i$ and the second case when $g$ starts with an
 element of $\Gamma^-_i$. Hence for all $g\in G$:
 $$\eqalign{\lambda(g_i)\delta_g &=
 \left\{\matrix{e^+_i\lambda(g_i)\delta_g
 &\hbox{if}&|g_ig| = |g|+1\cr&\cr
 \lambda(g_i) e^-_i\delta_g&\hbox{if}&|g_ig| =
|g|-1\cr}\right.\cr
 & = e^+_i\lambda(g_i)\delta_g + \lambda(g_i) e^-_i\delta_g \quad
 \hbox{(all
 cases).}}$$
 Therefore
 $$\lambda(g_i) = e^+_i\lambda(g_i) + \lambda(g_i)e^-_i$$
 and by taking adjoints:
 $$\lambda(g^{-1}_i) = e^-_i\lambda(g^{-1}_i) +
 \lambda(g^{-1}_i)e^+_i.$$
 Set
 $$\left.\eqalign{ u_i &= e^+_i\lambda(g_i), \quad u_{n+i} =
e^{-}_i
 \lambda(g^{-1}_i)\cr
 v_i &= \lambda(g^{-1}_i)e^-_i, \quad v_{n+i} =
 \lambda(g^{-1}_i)e^+_i}\right\} i=1,\ldots, n$$
 and for simplicity of notation, set also $g_{n+i} = g^{-1}_i$,
 $i=1,\ldots,
 n$. Then
 $$\lambda(g_i) = u_i+v_i,\qquad i=1,\ldots, 2n.$$
 Since $\sum\limits^n_{i=1} (e^+_i + e^-_i) = 1-e_0$ we have
 $$\sum^{2n}_{i=1} u_iu_i^* = \sum^{2n}_{i=1} v_i^*v_i =
 1-e_0 \le 1.$$
 So
 $$\left\|\sum^{2n}_{i=1} u_iu^*_i\right\| \le 1\quad \hbox{and}\quad
 \left\|\sum^{2n}_{i=1} v^*_iv_i\right\| \le 1.$$
 For elements $c_1,\ldots, c_m$, $d_1,\ldots, d_m$ of a $C^*$-algebra
 $B$
 one has easily that
 $$\left\|\sum^m_{i=1} c_id_i\right\| \le \left\|\sum^m_{i=1}
 c_ic^*_i\right\| ^{1/2} \left\|\sum^m_{i=1} d^*_id_i\right\|^{1/2}.$$
 Hence, with $u_1,\ldots,u_{2n}, v_1,\ldots, v_{2n}$ as above and
 $a_1,\ldots, a_{2n}\in A$,
 $$\eqalign{\left\|\sum^{2n}_{i=1} u_i\t a_i\right\|_{C^*_r(\F_n)
\t A} &= \left\|\sum^n_{i=1} (u_i\t
1) (1\t a_i)\right\|_{C^*_r(\F_n)
 \t  A}\cr
 &\le \left\|\sum
 u_iu^*_i\right\|^{1/2}\left\|\sum^{2n}_{i=1} a^*_ia_i\right\|^{1/2}
 \cr
 &\le \left\|\sum a^*_ia_i\right\|^{1/2}}$$
 and similarly
 $$\eqalign{\left\|\sum v_i\t a_i\right\|_{C^*_r(\F_n)\t A} &\le
 \left\|\sum^n_{i=1}  (1\t
a_i)(v_i\t 1)\right\|_{C^*_r(\F_n)\t A}\cr
 &\le  \left\| \sum
 a_ia^*_i\right\|^{1/2}\left\|\sum^n_{i=1} v^*_iv_i\right\|^{1/2}\cr
 &\le \left\|\sum a_ia^*_i\right\|^{1/2}}$$
 so altogether
 $$\eqalign{\left\|\sum^{2n}_{i=1} \lambda(g_i) \t a_i\right\| &=
 \left\|\sum^{2n}_{i=1} u_i\t a_i + \sum^{2n}_{i=1} v_i\t
 a_i\right\|\cr
 &\le \left\|\sum^{2n}_{i=1} a^*_ia_i\right\|^{1/2} + \left\|
 \sum^{2n}_{i=1} a_ia^*_i\right\|^{1/2}\cr
 &\le 2 \max\left\{\left\| \sum^{2n}_{i=1} a^*_ia_i\right\|^{1/2},
 \left\|\sum^{2n}_{i=1} a_ia^*_i\right\|^{1/2}\right\}.}$$
 This proves (2) in the case $S = \{g_1,\ldots, g_n, g^{-1}_1,\ldots,
 g^{-1}_n\}$, and the remaining cases follows from this by setting some
 of
 the $a_g$'s equal to 0.  \qed

 \n {\bf Remark.} The preceding statement remains true (with the
 obvious
 modifications) for the free group on infinitely
 many generators. See also Proposition 4.9 below for a
generalization of (1.1) and (1.2).

 \n {\bf Remark 1.2.} The proof of (2) is an illustration
of the following
 general principle. Let $T_1,\ldots, T_n$ be operators on a Hilbert
 space
 $H$ and let $c$ be a constant. The following properties
are essentially equivalent:\medskip
 \item{$(i)_c$} For any $C^*$-algebra $A$ and any set
$(a_i)_{i\le n}$ in $A$ we
 have
 $$\left\|\sum T_i\t a_i\right\| \le c\max
\left\{\left\|\left( \sum
 a^*_ia_i\right)^{1/2}\right\|, \left\|\left(\sum
 a_ia^*_i\right)^{1/2}\right\|\right\}.$$
 \item{$(ii)_c$} There are operators $u_i,v_i$ in $B(H)$
such that $T_i =
 u_i+v_i$ and
 $$\left\|\left(\sum u^*_i u_i\right)^{1/2}\right\| + \left\|\left(
 \sum
 v_iv^*_i\right)^{1/2}\right\| \le c.$$
 \medskip
 \n More precisely, we have $(ii)_c \Rightarrow (i)_c$ and
$(i)_c \Rightarrow (ii)_{2c}$
 \n The implication $(ii)_c \Rightarrow (i)_c$ follows as
above from the
 triangle inequality. To prove the converse, note that
$(i)_c$ equivalently means
 that the operator $u\colon\ E_n\to B(H)$ which maps $\delta_i$ to
 $T_i$
 satisfies $\|u\|_{cb}\le 1$.
 By the extension property of c.b. maps (cf.[Pa, p.100])
there
 is an extension $\tilde u\colon \ M_n\oplus M_n\to B(H)$ such that
 $\tilde
 u(\delta_i) = T_i$ and $\|\tilde u\|_{cb}\le 1$. Letting $u_i = \tilde
 u(e_{i1}\oplus 0)$ and $v_i = \tilde u(0\oplus e_{1i})$
 we obtain a
 decomposition satisfying $(ii)_{2c}.$ This shows that
$(i)_c$ implies $(ii)_{2c}$.

 \proclaim Proposition 1.3. Let $E_n \subset M_n \oplus M_n$ be the
 operator
 space
 $$E_n = \left\{\left(\matrix{c_1\cr \vdots&\bigcirc\cr c_n\cr}\right)
 \oplus \left(\matrix{c_1&\ldots&c_n\cr&\cr &\bigcirc\cr}\right)
 \ \Big| \
 c_1,\ldots, c_n \in {\bf C}\right\}.$$
 Then there are linear mappings
 $$w\colon \ E_n \to C^*_\lambda(\F_n) \quad \hbox{and}\quad v\colon \
 C^*_\lambda(\F_n) \to E_n$$
 such that
 $$vw = I_{E_n}\quad \hbox{and}\quad \|v\|_{cb} \|w\|_{cb}\le 2.$$
 Similarly, for the von~Neumann algebra $VN(\F_n)$ generated $\lambda$,
 there
 are linear mappings
 $$w_1\colon  \ E_n\to VN(F_n)\quad \hbox{and}\quad  v_1\colon
 \ VN(F_n) \to
 E_n$$
 such that
 $$v_1w_1 = I_{E_n} \quad \hbox{and}\quad \|v_1\|_{cb} \|w_1\|_{cb} \le
 2.$$
 In particular $E_n$ is $cb$-isomorphic to a $cb$-complemented subspace
 of
 $C^*_\lambda(\F_n)$ (resp. of $VN(\F_n)$).

 \n {\bf Proof.} Let $(\delta_1,\ldots, \delta_n)$ be the basis of
 $E_n$ determined by
 $$\sum^n_{i=1} c_i\delta_i = \left(\matrix{c_1\cr
\vdots&\bigcirc\cr
 c_n\cr}\right) \oplus \left(\matrix{c_1&\ldots&c_n\cr
 &\cr &\bigcirc\cr}\right)$$
 for $c_1,\ldots, c_n \in {\bf C}$. Define $w\colon \ E_n\to
 C^*_\lambda(\F_n)$ by
 $$w\left(\sum^n_{i=1} c_i \delta_i\right) = \sum^n_{i=1} c_i
 \lambda(g_i)$$
 and $v\colon \ C^*_\lambda(\F_n) \to E_n$ by
 $$v(x) = \sum^n_{i=1} \tau(\lambda(g_i)^* x)\delta_i$$
 where $\tau$ is the trace on $C^*_\lambda(\F_n)$ given by
 $$\tau(y) = (y\delta_e,\delta_e),\quad y\in C^*_\lambda(\F_n). \quad
 \hbox{(cf. [KR, p. 433])}.$$
 For any set $a_1,\ldots, a_n$ of $n$ elements in a $C^*$-algebra $A$
 $$(w\t I_A) \left(\sum^n_{i=1} \delta_i\t a_i\right) = \sum^n_{i=1}
 \lambda(g_i)\t a_i.$$
 Since
 $$\eqalign{\left\|\sum^n_{i=1} e_i\t a_i\right\| &=
 \left\|\left(\matrix{a_1\cr \vdots&\bigcirc\cr a_n\cr}\right) \oplus
 \left(\matrix{a_1&\ldots&a_n\cr&\cr &\bigcirc\cr}\right)\right\|\cr
 &= \max\left\{\left\|\sum^n_{i=1} a^*_ia_i\right\|^{1/2},
 \left\|\sum^n_{i=1} a_ia^*_i\right\|^{1/2}\right\}}$$
 it follows from Theorem 1.1 (2), that $\|w\t I_A\| \le 2$. Hence
 $\|w\|_{cb} \le 2$. Since
 $$\tau(\lambda(g)^* \lambda(h)) = \left\{\matrix{1&g=h\cr
0&g\ne
 h\cr}\right.$$
 we get for any finite subset $S\subset F_n$ and scalars $(c_g)_{g\in
 S}$
 $$v\left(\sum_{g\in S} c_g\lambda(g)\right) = \sum^n_{\scriptstyle i=1
 \atop \scriptstyle (g_i\in S)} c_{g_i}\delta_i$$
 and hence
 $$(v\t I_A) \left(\sum_{g\in S} \lambda(g) \t a(g)\right)
=
 \sum^n_{\scriptstyle i=1\atop \scriptstyle (g_i\in S)} \delta_i\t
 a(g_i).$$
 Let $S^1 = S\cap \{g_1,\ldots, g_n\}$. Then
 $$\eqalign{\left\|\sum^n_{\scriptstyle i=1\atop \scriptstyle g_i\in S}
 \delta_i\t a(g_i)\right\| &= \max\left\{\left\|\sum_{g\in S^1} a(g)^*
 a(g)\right\|^{1/2}, \left\|\sum_{g\in S^1} a(g) a(g)^*
 \right\|^{1/2}\right\}\cr
 &\le \max\left\{\left\|\sum_{g\in S} a(g)^* a(g)\right\|^{1/2},
 \left\|\sum_{g\in S} a(g) a(g)^*\right\|^{1/2}\right\},}$$
 which by Theorem 1.1(1) is smaller than or equal to
 $$\left\|\sum_{g\in S} \lambda(g) \t a(g)\right\| _{C^*_r(\F_n)
\t A}.$$
 Hence $\|v\t I_A\| \le 1$ and thus $\|v\| _{cb}\le 1$. Therefore
 $$\|v\|_{cb} \|w\|_{cb}\le 2$$
 and by construction $v  w = I_{E_n}$. This implies that
$w$ is a
 $cb$-isomorphism of $E_n$ onto its range
 $$w(E_n) = \hbox{span}\{\lambda(g_i) \mid i=1,\ldots, n\}$$
 and
 $$\|w\|_{cb} \|w^{-1}\|_{cb}\le 2.$$
 Moreover $P=wv$ is a completely bounded projection of
$C^*_\lambda(\F_n)$
 onto $w(E_n)$ and $\|P\|_{cb}\le 2$. The proof with
$VN(\F_n)$ in the place
 of $C^*_\lambda(\F_n)$ is easy since $v$ admits an
extension $v_1: VN(F_n)\to E_n$ with $\|v_1\|_{cb}
\le 1$. We leave the details
 to the reader. \qed

 \proclaim Lemma 1.4. ([P1, H1, LPP]). Let $u\colon \ A\to
B$ be a bounded
 linear operator between two $C^*$-algebras $A$ and $B$. Then for every
 $n
 \in {\bf N}$
 $$\|I_{E_n}\t u\|_{E_n\t A\to E_n\t B} \le \sqrt 2\, \|u\|.$$

 \n {\bf Proof.} The statement of the lemma is equivalent to:\ For all
 $a_1,\ldots, a_n \in A$
 $$ { \max\left\{\left\|\sum u(a_i)^*
u(a_i)\right\|, \left\| \sum
 u(a_i) u(a_i)^* \right\|\right\}
  \leq 2 \|u\|^2 \max\left\{\left\| \sum
a^*_ia_i\right\|, \left\|\sum
 a_ia^*_i\right\|\right\}.}\leqno(1.3) $$
 This is essentially [P1], (see also [H1,LPP]). However
to get the
 constant 2 in (1.3) one has to modify the proof of [H1,
Cor.~3.4] slightly:

 \n Let $T\colon \ A\to H$ be a bounded linear operator
from the $C^*$-algebra
 $A$ with values in a Hilbert space. By [H1, Thm.~3.2],
 $$\sum\|T(a_k)\|^2 \le \|T\|^2 \left(\left\| \sum
a^*_ka_k\right\| +
 \left\| \sum a_ka^*_k\right\|\right). \leqno (1.4)$$
 We can assume, that $B\subseteq B(K)$ for some Hilbert space $K$. By
 the
 above inequality (1.4) we get for any $\xi \in K$, that
 $$\sum \|u(a_k)\xi\|^2 \le \|\xi\|^2 \|u\|^2 \left(\left\| \sum
 a^*_ka_k\right\| + \left\| \sum a_ka^*_k\right\|\right).$$
 Clearly (1.4) also holds for conjugate linear maps, so
 $$\sum \|u(a_k)^*\xi\|^2 \le \|\xi\|^2 \|u\|^2 \left(\left\| \sum
 a^*_ka_k\right\| + \left\|\sum a_ka^*_k\right\|\right).$$
 Thus
 $$\max\left\{\left\| \sum u(a_k)^* u(a_k)\right\|, \left\|\sum u(a_k)
 u(a_k)^*\right\|\right\} \le \|u\|^2 \left(\left\| \sum
 a^*_ka_k\right\| +
 \left\|\sum a_ka^*_k\right\|\right)$$
 which implies (1.3).\qed

 \proclaim Theorem 1.5. Let $u\colon \ A\to B$ be a bounded linear
 operator
 between two $C^*$-algebras $A$ and $B$. Then for every
$k,n\in {\bf N}$
 $$\|I_{E^k_n}\t u\|_{E^k_n\t A\to E^k_n\t B} \le 2^{{3\over
 2}k-1}\|u\|.$$

 \n {\bf Proof.} The theorem is proved by induction on $k$.
By Lemma~1.4 the
 theorem holds for $k=1$. Assume next that the theorem is true for a
 particular $k\in {\bf N}$. Let
 $$w\colon \ E_n\to C^*_\lambda(\F_n)\quad \hbox{and}\quad v\colon \
 C^*_\lambda(\F_n) \to E_n$$
 be as in Proposition~1.2, and let $u\colon \ A\to B$ be a linear map
 between two $C^*$-algebras $A$ and $B$. Clearly
 $$I_{E_n}\t u = (v\t u) (w\t I_A)\leqno (1.5)$$
 where
 $$\eqalign{\|v\t u\| &= \|(v\t I_B)(I_{E_n}\t u)\|\cr
 &\le \|v\|_{cb} \|I_{E_n} \t u\|\cr
 &\le \sqrt 2\, \|u\| \|v\|_{cb}}$$
 by Lemma~1.4. Moreover $v\t u$ maps the $C^*$-algebra
 $C^*_\lambda(\F_n)\t
 A$ into the $C^*$-algebra $M_n(B) \oplus M_n(B)$, so by the induction
 hypothesis
 $$\|I_{E^k_n} \t v\t u\| \le 2^{{3\over 2}k-1} \|v\t u\| \le
 2^{{3\over
 2}k-{1\over 2}} \|u\| \|v\|_{cb}.$$
On the other hand by (0.5)
$$\|I_{E^k_n} \t w \t I_A \| = \| I_{{E^k_n}\t A} \t w \|
\le \| w\| _{cb}$$
 Now by (1.5)
 $$I_{E^{k+1}_n} \t u = (I_{E^k_n} \t v \t u)
(I_{E^k_n} \t w \t I_A).$$
 Thus, by Proposition~1.3
 $$\eqalign{\|I_{E^{k+1}_n} \t u\| &\le 2^{{3\over 2}k-{1\over 2}}
 \|u\|
 \|v\|_{cb} \|w\|_{cb}\cr
 &\le 2^{{3\over 2}k+{1\over 2}} \|u\|\cr
 &= 2^{{3\over 2}(k+1)-1} \|u\|.}$$
 Hence Theorem 1.5 follows by induction on $k$.\qed

\vfill\eject

\n {\bf \S 2. Description of $E_n^k$.}

In this section, we will identify the norm in the space $E_n^k
 {\otimes} A$  with the norm previously introduced in
(0.3) and (0.4) as $[\quad ]_{(k)}$.

\proclaim Proposition 2.1. Let $A$ be any $C^*$-algebra. Let
$n \geq 1,k \geq 1$ and let $\{a_J|J \in [n]^k\}$ be
elements of $A$. Then $$[(a_J)]_{(k)} = \Big|\Big| {\sum_{J
\in [n]^k} \delta_J \otimes a_J}\Big|\Big|_{E_n^k
  {\otimes} A} \leqno (2.1)$$
where we denote if $J = (j_1,...,j_k)$
$$\delta_J = \delta_{j_1} \otimes ... \otimes
\delta_{j_k}$$

\n The proof below is easy but the notation is a bit
painful. Using Proposition 2.1 we can complete the proof of
the results announced in the introduction.

\n {\bf Proof of Theorem 0.k} : Consider an operator $u : A
\to B$ between $C^*$-algebras. By Theorem 1.5 we have for
all $(a_J)$ in $A$ $$\N{\sum \delta_J \otimes
u(a_J)}_{E_n^k   {\otimes} B} \leq 2^{(3k/2)-1} \N {u} \N
{\sum \delta_J \otimes a_J}_{E_n^k   {\otimes} A}.$$
Taking (2.1) into account this immediately implies
$(0.1)_k$ and completes the proof of Theorem 0.k.

We   now check (2.1).
We will need the following elementary fact

\proclaim Lemma 2.2. Let $H,H_1,H_2,H_3,H_4$ be Hilbert
spaces.  Let $e \in H_1,f \in H_4$ be norm one vectors.
Let $(\ph_j)_{j \in J}$ and $(\psi_i)_{i \in I}$ be orthonormal finite
sequences in $H_2$ and $H_3$ respectively. Let $a_{ij}$ be elements of
a
$C^*$-algebra $A$ embedded into $B(H)$.
Then we have
$$\Big|\Big| {\sum_{\buildrel{i \in I}\over {j \in J}}(e \otimes
\ph_j) \otimes (\psi_i \otimes f) \otimes a_{ij}}\Big|\Big|
= \sup_{\buildrel{y_i \in H}\over {x_j \in H}}
\Big{\{}\Big|\sum_{i,j} < y_i,a_{ij}-x_j> \Big| \ ,\  \sum
\| {x_j}\|^2 \leq 1,\sum \| {y_i}\|^2 \leq
1\Big\}.\leqno(2.2) $$ Here the norm on the left hand side
means the norm in the space of all bounded operators from
$H_1 \otimes_2 H_2 \otimes_2 H$ into $H_3 \otimes_2 H_4
\otimes_2 H$.

\n {\bf Proof}. We may clearly assume without loss of
generality that $H_1 = \comp e,H_4 = \comp f$ and that
$(\ph_j)$ (resp. $(\psi_i))$ is a basis of $H_2$ (resp.
$H_3$).
Then the norm we want to compute is clearly equal to the norm of the
operator
$$\tilde T = \sum_{ij} \ph_j \otimes \psi_i \otimes a_{ij}$$
as an operator from $H_2 \otimes_2 H$ to $H_3 \otimes_2 H$. But then
the
general form of an element in the unit ball of $H_2 \otimes_2 H$ (resp.
$H_3
\otimes_2 H)$ is given by $\sum \ph_j \otimes x_j$ (resp. $\sum \psi_i
\otimes
y_i)$ with $x_j \in H_2$ (resp. $y_i \in H_3)$ such that $\sum
\| {x_j}\|^2 \leq 1$ (resp. $\sum \| {y_i}\|^2 \leq 1)$.
Hence the norm of $\tilde T$ (or of $T$) is equal to the
right hand side of (2.2).\qed

We need to introduce more notation.

\n Recall that $E_n \subset M_n \oplus M_n$ and $\delta_i =
e_{i1}-\oplus e_{1i}$. We consider of course $M_n \oplus
M_n$ as a subset of the set of all operators on $\ell_2^n
\oplus \ell_2^n$.  It will be convenient to denote
$e_{ij}^0 = e_{ij} \oplus 0$ and $e_{ij}^1 = 0 \oplus
e_{ij}$ in $M_n \oplus M_n$. Also
$e_i^0 = e_i \oplus 0$ and $e_i^1 = 0 \oplus e_i$ in $\ell_2^n \oplus
\ell_2^n$. As is usual, for $e,f$ in $H$,  we will identify the tensor
$e \otimes
f ~- $ with the operator $x \to < e,x> f$ ( defined on $ H)$.
Hence in tensor product notation we have (with the usual matricial
conventions) $e_{ij} = e_j \otimes e_i$ and $\delta_i = e_1^0 \otimes
e_i^0 +
e_i^1 \otimes e_1^1$.
Let us denote by $H_0$ the span of $\{e_{i1}| i = 1,...,n\}$ in $M_n$
and by
$H_1$ the span of $\{e_{1i}|i=1,...,n\}$ in $M_n$, so that $E_n \subset
H_0
\oplus H_1$.
Let $P_0 : H_0 \oplus H_1 \to H_0$ (resp. $P_1 : H_0 \oplus H_1 \to
H_1)$
denote the canonical projection.
We have $E_n^{\otimes k} \subset (H_0 \oplus H_1)^{\otimes k}$.
For $\alpha \in \{0,1\}^k$ we denote
$$P_\alpha : (H_0 \oplus H_1)^k \to (H_0 \oplus H_1)^k$$
the projection defined by
$$P_\alpha = P_{\alpha(1)} \otimes P_{\alpha(2)} \otimes...\otimes
P_{\alpha(k)}.$$
Let us denote by $I_X$ the identity on $X$.
Then we have
$$\eqalign{
I_{(H_0\oplus H_1)^{\otimes k}} & = (I_{{H_0}\oplus H_1})^{\otimes
k}\cr
& = (P_0 + P_1)^{\otimes k}\cr
& = \sum_{\alpha \in \{0,1\}^k} P_{\alpha(0)}\otimes...\otimes
P_{\alpha(k)}\cr
& = \sum_\alpha P_\alpha}\leqno (2.3)$$

\n {\bf Proof of Proposition 2.1} : Let $T = \sum_{J \in
{[n]}^k} \delta_J \otimes a_J$. By (2.3) we have
$$T = \sum_\alpha T_\alpha$$
where
$$T_\alpha = \sum_J P_\alpha (\delta_J) \otimes a_J.$$
We now claim that
$$\| {T_\alpha}\| = \| {(a_J)}\|_\alpha.
 \leqno (2.4)$$
To check this, we can assume for simplicity (up to a permutation of the
factors
in the tensor product) that $\alpha$ is the indicator function of the
set
$\{1,2,...,p\}$ for some $p$ with $1 \leq p \leq k$. Then if $J =
(j_1,...,j_k)$ we have
$$P_\alpha (\delta_J) = e_{j_1 1}^1 \otimes...\otimes e_{j_p1}^1
\otimes
e_{1j_{p+1}}^0 \otimes...\otimes e_{1j_k}^0.\leqno(2.5)$$
(Recall the convention that the tensor $e \otimes f$ represents the
operator $x
\to < e,x> f)$.
Let $e^1(\alpha) = e_1^1 \otimes...\otimes e_1^1$~~~~~~~~~-($p$
times) and $f^0(\alpha) = e_1^0 \otimes...\otimes
e_1^0$~~~~~~~~($k-p$ times).

\n Then (2.5) yields
$$P_\alpha (\delta_J) = (e^1(\alpha) \otimes e^0_{j_{p+1}}
\otimes...\otimes
e^0_{j_k}) \otimes (e^1_{j_1} \otimes...\otimes e^1_{j_p} \otimes
f^0(\alpha)).$$
If we now write $e^\eps_{\{j_1,...,j_p\}}$ instead of
$e^\eps_{j_1}\otimes...\otimes e^\eps_{j_p}$ for $\eps = 0$ or $1$, we
can
rewrite the last identity as
$$P_\alpha(\delta_J) = (e^1(\alpha) \otimes e^0_{{\pi_\alpha}(J)})
\otimes
(e^1_{{\pi_{\alpha^c}}(J)} \otimes f^0(\alpha)),\leqno(2.6)$$
where we recall that $\pi_\alpha : [n]^k \to [n]^\alpha$ denotes the
canonical
projection.
Then the above lemma 2.2 gives in the present particular case
$$ {
\| {T_\alpha}\|   = \big{|}\big{|} {\sum_J (e^1(\alpha)
\otimes e^0_{{\pi_\alpha}(J)})\otimes
(e^1_{{\pi_{\alpha^c}}(J)}\otimes f^0(\alpha))\otimes
a_J} \big{|}\big{|} = \| {(a_J)}\|_\alpha}.$$
 This proves our
claim (2.4).

\n Now, we can conclude.
Let us denote $h^0 = \ell_2^n \oplus 0$ and $h^1 = 0 \oplus \ell_2^n$
in
$\ell_2^n \oplus \ell_2^n$. Let $K_\alpha$ be the support of $T_\alpha$
(\ie
the orthogonal of its kernel) and let $R_\alpha$ be the range of
$T_\alpha$.
Then the preceding formula (2.6) shows that $K_\alpha$ is equal to the
tensor
product $F_1 \otimes F_2 \otimes...\otimes F_k$
where
$$F_j = \comp e^1_1~- {\rm if}~-j \in \alpha$$ and
$$F_j = h^0~-{\rm if}~-j \not\in \alpha.$$
It follows that the subspaces $(K_\alpha)$ are mutually orthogonal.
Similarly,
the family $(R_\alpha)$ is mutually orthogonal. By a well known
estimate it
follows that
$$\N {\sum T_\alpha} = \max_\alpha \N {T_\alpha}.$$
This completes the proof. \qed

\vfill\eject

 \overfullrule = 0pt

 \n {\bf \S 3. Random series in non-commutative $L_1$-spaces.}

 Let $A$ be a von Neumann algebra with a predual denoted by $A_*$.

 \n Let $\xi_1,\ldots, \xi_n \in A_*$ and let (recall
(0.2))
 $$[(\xi_i)]^*_{(1)} = \sup \left\{\left| \sum \langle \xi_i,a_i\rangle
 \right| \ \bigg| \ a_i \in A \ \ [(a_i)]_{(1)} \le
1\right\}.$$
 For instance, if $A=B(H)$, $A_* = C_1(H)$ (the space of trace class
 operators on $H$) and we have clearly
 $$[(\xi_i)]^*_{(1)} = \inf\left\{ {tr} \left(\sum
 x^*_ix_i\right)^{1/2} +
 {tr}\left(\sum y_iy^*_i\right)^{1/2}\right\}$$
 where the infimum runs over all decompositions $\xi_i = x_i +y_i$ in
 $C_1(H)$.

 Let ${\bf T}^{\bf N}$ be the  infinite dimensional torus equipped with
 its
 normalized Haar measure $\mu$. The following result is proved in
 [LPP].

 \n For all $\xi_1,\ldots, \xi_n$ in $A_*$
 $${1\over 2} [(\xi_i)]^*_{(1)} \le \int \left\|\sum^n_{j=1}
 e^{it_j}\xi_j
 \right\|_{A_*} d\mu(t) \le [(\xi_i)]^*_{(1)}.\leqno (3.1)$$
 (See Theorem 3.3 below and its proof.)

 It is easy to deduce from (3.1) a necessary and sufficient condition
 for a
 series of the form
 $$S(t) = \sum^\infty_{j=1} e^{it_j} \xi_j\ ,\qquad
t=(t_j)_{j\in {\bf N}} \in {\bf T}^{\bf N}$$
 to converge in $L_2({\bf T}^{\bf N}, \mu; A_*)$. The aim
of this section is
 to prove a natural extension   of (3.1) to double
 series of the form
 $$S(t',t'') = \sum^\infty_{j,k=1} e^{it'_j} e^{it''_k} \xi_{jk}$$
 with $\xi_{jk}\in A_*$, $t',t'' \in {\bf T}^{\bf N}$. More generally,
 we
 will consider for any $k\ge 1$, elements $\xi_{j_1j_2\ldots j_k}$ in
 $A_*$
 and  will find an equivalent for the expression
 $$\int \left\|\sum_{j_1\le n,\ldots, j_k\le n} e^{it^1_{j_1}}\ldots
 e^{it^k_{j_k}} \xi_{j_1j_2\ldots j_k} \right\|_{A_*} d\mu(t^1) \ldots
 d\mu(t^k).$$
 See Theorem 3.6 below for an explicit statement.

 Let $A$ be a $C^*$-algebra throughout this section. We will denote
 simply
 $$\eqalignno{C_n &= C^*_\lambda(\F_n)\cr
 \noalign{\hbox{and}}
 C^k_n &= C_n\t \cdots \t C_n \qquad (k \hbox{ times}).}$$
 We always equip the tensor products such as $E_n\t A$, $C_n\t A$,
 $C^k_n\t
 A$ with the spatial (or minimal) tensor product. More precisely,
 whenever
 $S\subset B(K)$ is an operator space and $A\subset B(H)$ is a
 $C^*$-algebra, we will denote by $S\t A$ the linear tensor product
 equipped
 with the norm induced by $B(K \t_2H)$.

 Let $G$ be a discrete group. For $t\in G$, let $\lambda_*(t)$ denote
 the element of $C^*_\lambda(G)^*$ given by
 $$\forall\ a\in C^*_\lambda(G)\qquad \langle \lambda_*(t), a\rangle =
 \langle
 a\delta_e,\delta_t\rangle.$$
 Clearly $$\langle \lambda_*(s), \lambda(t)\rangle=
 \left\{\matrix{1&\hbox{if } s=t\hfill\cr&\cr
 0&\hbox{otherwise.}\hfill\cr}\right.$$

\n Note that if
$C^*_\lambda(G)^*$ is
 identified with $B_\lambda(G)$ in the usual way (see for instance
[E])
 then $\lambda_*(t)$ simply corresponds to the function $\delta_t$.

 \n For any $J = (j_1,\ldots, j_k) \in [n]^k$ we denote by
$g_J$ the element
 of $(\F_n)^k$ defined by
 $$g_J = (g_{j_1},\ldots, g_{j_k}).$$
 Then with the obvious identification
 $$C^*_\lambda((\F_n)^k) = C^k_n$$
 we have $\lambda(g_J) = \lambda(g_{j_1})\t \cdots \t
 \lambda(g_{j_k})$. We
 will also consider the dual $E^*_n$ of the space $E_n$ considered in
 section~1 and  will denote by $\{\delta^*_j\}$ the basis
of $E^*_n$ which
 is biorthogonal to $\{\delta_j\}$.
 We will also consider $E^k_n = E_n\t \cdots \t E_n$ ($k$ times) and
 its
 dual $(E^k_n)^*$. We will denote for any $J = (j_1,\ldots, j_k)$ in
 $[n]^k$
 $$\delta^*_J = \delta^*_{j_1}\t \cdots \t
\delta^*_{j_k}\in (E^k_n)^*$$ and
$$\lambda_*(g_J) = \lambda_*(g_{j_1})\t \cdots \t
\lambda_*(g_{j_k})\in (C^k_n)^*.$$
 We will denote by $\Omega$ the infinite dimensional torus i.e. we set
 $$\Omega = {\bf T}^{\bf N}$$
 and we equip $\Omega$ with the normalized Haar measure
$\mu$.
 (In most of what follows, it would be more appropriate to replace
 $\Omega$
 by $\Omega_n = {\bf T}^n$, but we try to simplify the notation.)
 We will denote by
 $$\varepsilon_j \colon \ \Omega \to {\bf T}$$
 the sequence of the coordinate functions on $\Omega$. Moreover, we
 will
 consider the product $\Omega^k$ equipped with the product measure
 $\mu^k$.
 For any $J = (j_1,\ldots, j_k) \in [n]^k$, let $\varepsilon _J\colon \
 \Omega^k\to {\bf T}$ be the function defined by
 $$\forall\ (t_1,\ldots, t_k) \in \Omega^k \qquad
 \varepsilon_J(t_1,\ldots,
 t_k) = \varepsilon_{j_1}(t_1)\ldots \varepsilon_{j_k}(t_k).$$
 Equivalently $\varepsilon_J = \varepsilon_{j_1} \t \cdots \t
 \varepsilon_{j_k}$. We first record a simple consequence of
 Proposition~1.3.

 \proclaim Lemma 3.1. For any $\{\varepsilon_j \mid j\le n\}$ in $A^*$
 we
 have
 $${1\over 2} \left\|\sum \delta^*_j\t \xi_j\right\|_{(E_n\t A)^*} \le
 \left\|\sum \lambda_*(g_j)\t \xi_j\right\| _{(C_n\t A)^*}
\le \left\|\sum
 \delta^*_j \t \xi_j\right\|_{(E_n\t A)^*}.$$

 \n {\bf Proof.} Let $v,w$ be as in Proposition 1.3. Since
$w \delta_j =
 \lambda(g_j)$ and $v(\lambda(g_j)) = \delta_j$ we have
$(w\t I_A)^*
 (\lambda_*(g_j)\t \xi) = \delta^*_j\t \xi_j$ and $(v\t
I_A)^* (\delta^*_j\t
 \xi_j) = \lambda_*(g_j) \t \xi_j$.
 Hence, recalling (0.5),  Lemma~3.1 follows from
$\|w\|_{cb}\le 2$ and $\|v\|_{cb}\le
 1$.\qed

 The next lemma is rather elementary.

 \proclaim Lemma 3.2.  \item {\rm (i)} Consider $\{\xi_{ij}\mid i,
j=1,...,n  \}$ in $A^*$. For any orthonormal systems
$\varphi_1,...,\varphi_n$ and $\psi_1,...,\psi_n$ in
$L_2(\mu)$ (where $\mu$ is a probability as above) we have
$$\int \| \sum \varphi_i(t)\psi_j(s)
\xi_{ij}\|_{A^*}d\mu(t)d\mu(s) \leq
\|(\xi_{ij})\|_{M_n(A)^*}.\leqno(3.2)$$
 \item {\rm (ii)} For any $k\ge 1$ and any $(\xi_J)$ in
$A^*$ we have
$$\left\|\sum_{J\in [n]^k}
\varepsilon_J\xi_J\right\|_{L_1(\mu;A^*)}\leq
\left\|\sum_{J\in [n]^k}
 \delta^*_J \t \xi_J\right\|_{(E_n^k\t A)^*}.\leqno(3.3)$$

\p (i) To prove this, it clearly suffices to assume that
$A$ is a von Neumann algebra and that $\xi_{ij}\in A_*$.
Since $M_n(A)$ is a subspace of $M_n(B(H))$ for some
Hilbert space $H$, by duality   its predual  $M_n(A)_*$
is a quotient of $M_n(B(H))_*$. This shows that it
suffices to prove (i) for $A=B(H)$ and $\xi_{ij}\in
B(H)_*$. Then we can identify $M_n(B(H))_*$ with the
projective tensor product $\ell_2^n(H)^*\hat{\t}
\ell_2^n(H)$. Consider an element $x$ (resp. $y$) in the
unit ball of
$\ell_2^n(H)$ (resp. $\ell_2^n(H)^*$). Let $\xi$ be the
element of $M_n(B(H))_*$ defined by $\xi=y\t x$ or
equivalently, $\xi= (\xi_{ij})$ with $\xi_{ij}=y_j\t x_i$.
For such a $\xi$ we have
$$(\int \| \sum   \varphi_i(t)\psi_j(s)
\xi_{ij}\|^2_{A^*}d\mu(t)d\mu(s))^{1/2} =
(\int \| \sum \varphi_i(t) x_i\|^2 d\mu(t)
\int \|\sum  \psi_j(s) y_j\|^2 d\mu(s))^{1/2}$$
$$
\ \ \ \ \ \ \ =\|x\|\  \|y\| \leq 1.$$
Since the unit ball of $M_n(B(H))_*$  is the closed
convex hull of elements of this form, we obtain (3.2).

\n (ii) Consider a subset $\alpha\subset \{1,...,k\}$. We
denote by $\alpha^c$ its complement. Recall that for
elements  $(a_J)_{J\in [n]^k}$ in $A$ the norm
$\|(a_J)\|_\alpha$ defined in (0.3) can be viewed as the
norm of a matrix acting from $\ell_2([n]^{\alpha},H)$ into
$\ell_2([n]^{{\alpha}^c},H)$. Therefore we deduce from
(3.2) that for any $(\xi_J)_{J\in [n]^k}$ in $A^*$ we have
$$\int \left\|\sum_{J\in [n]^k}
\varepsilon_J\xi_J \right\|_{A^*}d\mu^k\leq \|
(\xi_J)\|_{\alpha}^*. \leqno(3.4)$$
Observe that by duality (2.1) has the following
consequence.

\n  If $\left\|\sum_{J\in [n]^k}
 \delta^*_J \t \xi_J\right\|_{(E_n^k\t A)^*}\leq 1$ then
there is a decomposition $$\xi_J=\sum_{\alpha \subset
\{1,...,k\}} \xi_J^\alpha\quad  {\rm  with}\quad
\sum_{\alpha}\|(\xi_J^\alpha)\|_\alpha^*\leq 1.$$Therefore
(3.3) follows from (3.4) and the triangle inequality.
\qed

 We now reformulate the main result of [LPP] in our framework.

 \proclaim Theorem 3.3. For any $\{\xi_j\mid j\le n\}$ in $A^*$ we have
 $$\left\|\sum \varepsilon _j\xi_j\right\|_{L_1(\mu; A^*)} \le
 \left\|\sum
 \delta^*_j \t \xi_j\right\|_{(E_n\t A)^*} \le
2\left\|\sum \varepsilon_j
 \xi_j\right\|_{L_1(\mu; A^*)}.\leqno (3.5)$$

 \n {\bf Proof.} The left side is
(3.3)
 above for $k=1$. By our earlier analysis of $E_n\t A$, the
right side is clearly
 equivalent to the following fact.

 \n Assume $\big\| \sum
\varepsilon_j\xi_j\big\|_{L_1(\mu;A^*)} <1$. Then there
 is a decomposition $\xi_j = x_j+y_j$ in $A^*$ such that
 $$\eqalign{\forall (a_j)\in A\qquad \left|\sum \langle
x_j,a_j\rangle\right| &\le
 \left\|\left(\sum a^*_ja_j\right)^{1/2}\right\| \cr
{\rm and}\ -\ \ \qquad
 \left|\sum \langle y_j,a_j\rangle \right| &\le
\left\|\left( \sum
 a_ja^*_j\right)^{1/2}\right\|. }$$
 This is precisely what is proved in section~II of [LPP], except that
 the
 sequence $(\varepsilon_j)$ on $\Omega$ is replaced by the sequence
 $(e^{i3^jt})$ on the one dimensional torus. By a routine averaging
 argument, one can then obtain the preceding fact as stated above with
 $(\varepsilon_j)$. (Note actually that the approach of [LPP] can be
 developed directly for the functions $(\varepsilon_j)$, this is
explicited
 in [P2].)

 We now relate certain series on ${\bf Z}^n$ (formed by iterating the
 expressions appearing in Theorem~3.3) with the corresponding series on
 the
 free group $\F_n = {\bf Z} *\cdots * {\bf Z}$. In other words, our aim
 is
 to compare for these series the free group $\F_n$ with $n$ generators
 with
 its commutative counterpart ${\bf Z}^n$.

 \proclaim Lemma 3.4. For any $\{\xi_J\mid J\in [n]^k\}$ in $A^*$ we
 have
 (the summation being over all $J$ in $[n]^k$)
 $$2^{-k} \left\|\sum
\varepsilon_J\xi_J\right\|_{L_1(\mu^k;A^*)} \le\left\| \sum
 \lambda_*(g_J)\t \xi_J\right\|_{(C^k_n\t A)^*} \le 2^k\left\|\sum
 \varepsilon
 _J\xi_J\right\|_{L_1(\mu^k;A^*)}.$$

 \n {\bf Proof.} By the preceding three statements, we know that this
 holds
 for $k=1$. We now argue by induction. Assume Lemma~3.4 proved for an
 integer $k\ge 1$, and let us prove it for $k+1$.

 \n Consider elements $\{\xi_{Jj} \mid J \in [n]^k, j\in
[n]\}$ in $A^*$. We
 have
 $$\sum_{J'\in [n]^{k+1}} \lambda_*(g_{J'}) \t \xi_{J'} = \sum_{J\in
 [n]^k}
 \lambda_*(g_J) \t \left(\sum_{j\le n} \lambda_*(g_j) \t
 \xi_{Jj}\right).$$
 By the induction hypothesis, we have
 $$\left\|\sum_{J'} \lambda_*(g_{J'}) \t \xi_{J'}\right\|_{(C^{k+1}_N\t
 A)^*}
 \le 2^k \int_{\Omega^k} \left\|\sum \varepsilon _J(t) \eta_J\right\|
 _{(C\t
 A)^*} d\mu^k(t)\leqno (3.6)$$
 where $\eta_J = \sum\limits_j \lambda_*(g_j)\t \xi_{Jj}$.

 Now for each fixed $t$ in $\Omega^k$, we have by (3.5) and Lemma~3.1
 $$\left\|\sum \varepsilon_J(t)\eta_J\right\|_{(C\t A)^*} \le 2\int
 \left\|\sum \varepsilon_J(t) \left(\sum \varepsilon_j(s)
 \xi_{Jj}\right)\right\|_{A^*} d\mu(s).$$
 Integrating over $t\in \Omega^k$ this yields
 $$\int\left\|\sum \varepsilon_J(t)\eta_J\right\|_{(C\t A)^*} d\mu^k(t)
 \le
 2 \int \left\|\sum_{J'\in [n]^{k+1}} \varepsilon_{J'}
 \xi_{J'}\right\|_{A^*} d\mu^{(k+1)},\leqno (3.7)$$
 hence (3.6) and (3.7) yield the induction step for $k+1$. This
concludes
 the proof for the right side inequality in Lemma
3.4. The proof of the other inequality is entirely
similar.\qed

 We now  come to the main result of this section.

 \proclaim Theorem 3.5. For any $\{\xi_J\mid J\in [n]^k\}$ in $A^*$, we
 have
 $$\left\|\sum \varepsilon_J\xi_J\right\|_{L_1(\mu^k,A^*)} \le
 \left\|\sum
 \delta^*_J \t \xi_J\right\|_{(E^k_n\t A)^*} \le 2^{2k}
\left\|\sum
 \varepsilon_J\xi_J\right\|_{L_1(\mu^k,A^*)}.$$

 \n {\bf Proof.} With $v$ and $w$ as in Proposition 1.1, we have $\|w
 ^{\t
 k}\|_{cb} \le 2^k$, hence by (0.5)
 $$\|w^{ \t k} \t I_A\|_{E^k_n\t A \to C^k_n \t A} \le
2^k.$$
 Moreover, we have $w^{\t k} (\delta_J) = \lambda(g_J)$ hence
$(w^{\t k} \t
 I_A)^* (\lambda_*(g_J) \t \xi_J) = \delta^*_J \t \xi_J$.
This yields
 $$\left\|\sum \delta^*_J\t \xi_J\right\|_{(E^k_n\t A)^*} \le 2^k
 \left\|\sum \lambda_*(g_J)\t \xi_J\right\|_{(C^k_n\t
A)^*}.$$
 Combined with Lemma 3.4 this gives the right side in Theorem 3.5. The
 left
 side  has already been proved in Lemma 3.2.
 \qed

 \n {\bf Remark.} A slight modification of our proof yields Theorem~3.5
 with
 the constant $2^{2k-1}$ instead of $2^{2k}$.

\n{\bf Remark.} Let $k$ be a fixed integer. Consider
the mapping
$$Q_k\colon C(\Omega^k)\to E_n^k$$
defined by
$$\forall f\in C(\Omega^k)\quad Q_k(f)=\sum_{J\in [n]^k}
\hat{f}(J) \delta_J,$$
where $\hat{f}$ is the Fourier transform of $f$, i. e.
$\hat{f}(J)=\int f(t) \bar{\epsilon}_J(t) d\mu^k(t).$
Let $N_k=Ker(Q_k)$. Dualizing (3.3) we find that
$\|Q_k\|_{cb}\le 1$. Hence, considering $Q_k$ modulo its
kernel and equipping $C(\Omega^k)/N_k$ with its quotient operator space
structure (in the sense of [BP,ER]), we find a
map  $$U_k: C(\Omega^k)/N_k\to E_n^k\quad {\rm with}\quad
\|U_k\|_{cb}\le 1.$$ Then  Theorem 3.5 admits the
following dual reformulation: $U_k: C(\Omega^k)/N_k\to
E_n^k$ is a complete isomorphism and $\|U_k^{-1}\|_{cb}\le
2^{2k}$. In other words,
the space
$C(\Omega^k)/N_k$  is, for each $k$, completely
isomorphic (uniformly with respect to $n$) to $ E_n^k$.

 Assume now that $A$ is a von~Neumann algebra and let $A_*$ be its
 predual.
 We define for any family $(x_J)_{J\in [n]^k}$ in $A_*$ the norm which
 is
 dual to the norm $\|\ \|_\alpha$ defined in (0.3). We set
 $$\|(x_J)\|^*_\alpha = \sup\left\{\left| \sum_{J\in
[n]^k} \langle a_J,
 x_J\rangle  \right|\quad a_J\in A, \quad \|(a_J)\|_\alpha
\le 1\right\}.\leqno (3.8)$$
 Then we define
 $$[(x_J)]^*_{(k)} =  \inf \sum_{\alpha \in \{0,1\}^k}
 \|x^\alpha_J\|^*_\alpha \leqno (3.9)$$
 where the infimum runs over all $x^\alpha_J$ in $A_*$ such that $x_J =
 \sum\limits_{\alpha \in \{0,1\}^k} x^\alpha_J$.

 \n Assume that $A =  (A_*)^*$ is a von~Neumann subalgebra
of $B(H)$ and let
 $q\colon \ N(H) \to A_*$ be the quotient mapping which is the
 preadjoint of
 the embedding $A\hookrightarrow B(H)$.  We can also write
 $$\|(x_J)\|^*_\alpha = \inf \left\{\sum |\lambda_m|\right\}\leqno
 (3.10)$$
 where the infimum runs over all the possibilities to write $(x_J)$ as
 a
 series
 $$x_J = \sum_m \lambda_m h^m_{\pi_\alpha(J)} \t
 k^m_{\pi_{\alpha^c}(J)}$$
 where $(h^m_i)_{i\in [n]^\alpha}$ and $(k^m_j)_{j \in [n]^{\alpha^c}}$
 are
 elements of $H$ such that $\sum\limits_i \|h^m_i\|^2 \le 1$ and
 $\sum\limits_j \|k^m_j\|^2 \le 1$ for each $m$.

 \n The identity of (3.8) and (3.10) is clear since the
dual norms are the same
 by (0.3). Similarly it is clear that the dual space to
$(A_*)^{n^k}$
 equipped with the norm $[ \quad]^*_{(k)}$ can be identified with
 $(A)^{n^k}$ equipped with the norm $[\quad]_{(k)}$. By
 Proposition~2.1,
 this means that $(A_*)^{n^k}$ equipped with the norm $[\quad]^*_{(k)}$
 can
 be viewed as a predual (isometrically) of $E^k_n\t A$. Hence, we can
 now
 rewrite Theorem 3.5 a bit more explicitly. For all $(x_J)$
in $A_*$, we have (as
 announced in the beginning of this section)
 $$(2 ^{2k})^{-1} [(x_J)]^*_{(k)} \le \left\|\sum
\varepsilon_J
 x_J\right\|_{L_1(\Omega^k, A_*)} \le [(x_J)]^*_{(k)}. \leqno(3.11)$$
 In particular, we can state  for emphasis.

 \proclaim Theorem 3.6. Let $A\subset B(H)$ be a von~Neumann subalgebra
 with
 predual $A_*$ and let $q\colon \ H\hat{\t} H\to A_*$ be
the corresponding
 quotient mapping. Consider $\{x_J\mid J\in [n]^k\}$ in $A_*$ such that
 $$\left\|\sum \varepsilon_J x_J\right\|_{L_1(\Omega^k,A_*)} <1.$$
 Then $(x_J)$ admits a decomposition as
 $$x_J = \sum _{\alpha \in \{0,1\}^k} x^\alpha_J$$
 with
 $$x^\alpha_J = q\left(\sum_n \lambda^\alpha_m h^m_{\pi_\alpha(J)} \t
 k^m_{\pi_{\alpha^c}(J)}\right)$$
 where for each $\alpha$,  $\quad \{h^m_i \mid i \in
[n]^\alpha\}$ and $\{k^m_j \mid
 j\in [n]^{\alpha^c}\}$ are elements of $H$ such that $\sum\limits_i
 \|h^m_i\|^2 \le 1$ and $\sum\limits_j \|k^m_j\|^2 \le 1$ and where
 $\lambda^\alpha_m$ are scalars such that
 $$\sum_\alpha \sum_m |\lambda^\alpha_m| < 2^{2k}.$$
 Conversely, if $(x_J)$ admits such a decomposition, we
must have
 $\big\|\sum \varepsilon _Jx_J\big\|_{L_1(\Omega^k,A_*)} <
2^{2k}$.

 \n {\bf Proof.} The proof is nothing but (3.9), (3.10) and (3.11)
 spelt out
 explicitly.

 \n {\bf Remark.} The preceding theorem proves one of the conjectures
 formulated in [P2] in the case $A = B(H)$, $A_* = H\hat\t
H$.

\vfill\eject

 \n {\bf \S 4. Complements.}

 The following result shows that in Proposition~1.3, the
algebra
 $(C^*_\lambda(\F_n))^\infty_{n=1}$ cannot be substituted by any
 sequence of
 nuclear algebras.

 \proclaim Theorem 4.1. Let $A$ be either a nuclear $C^*$-algebra or an
 injective von~Neumann algebra, and let $I_{E_n} = vw$ be a
 factorization of
 $I_{E_n}$ through $A$. Then
 $$\|v\|_{cb} \|w\|_{cb} \ge {1\over 2}(1+\sqrt n).$$

 For the proof we need the following.

 \proclaim Lemma 4.2. Consider the subspace
 $S _n$ of $M_n\oplus M_n$ given
 by
 $$S_n = \left\{\left(\matrix{x_1\cr  x_2\cr
\vdots&\bigcirc\cr
 x_n\cr}\right) \oplus \left(\matrix{y_1&\ldots&y_n\cr
 &\cr &\bigcirc\cr}\right)
 \ \bigg| \ x_1,\ldots, x_n, y_1,\ldots, y_n \in {\bf C}\right\}$$
 and define $R\colon \ S_n\to S_n$ by
 $$R(x\oplus y) = y^t \oplus x^t,\qquad x\oplus y \in
S_n.$$
 Then\medskip
 \item{(a)} ${1\over 2}(I_{S_n}+R)$ is a projection of
$S_n$ onto $E_n$
 and
 $$\left\|{1\over 2} (I_{S_n}+R)\right\|_{cb} = {1\over
2}(1+\sqrt n).$$
 \item{(b)} For any projection $Q$ of $S_n$ onto $E_n$
(resp.
 $M_n\oplus M_n$ onto $E_n$) one has
 \medskip
 $$\|Q\|_{cb} \ge {1\over 2}(1+\sqrt n).$$

 \n {\bf Proof.} a)\ Obviously $R^2 = I_{ S_n }$ and $E_n = \{a\in  S_n
 \mid
 Ra=a\}$. Hence ${1\over 2}(I_{ S_n }+R)$ is a projection of $ S_n $
 onto
 $E_n$. Let $A$ be a $C^*$-algebra. Then
 $$ S_n  \t A = \left\{\left(\matrix{a_1\cr \vdots&\bigcirc\cr
 a_n\cr}\right)
 \oplus \left(\matrix{b_1&\ldots&b_n\cr &\cr&\bigcirc\cr}\right)\
 \bigg| \
 a_1,\ldots, a_n, b_1,\ldots, b_n \in A\right\}$$
 and
 $$(R\t I_A) \left(\left(\matrix{a_1\cr \vdots&\bigcirc\cr
 a_n\cr}\right)
 \oplus \left(\matrix{b_1&\ldots&b_n\cr &\cr
 &\bigcirc\cr}\right)\right) =
 \left(\matrix{b_1\cr\vdots&\bigcirc\cr b_n\cr}\right) \oplus
 \left(\matrix{a_1&\ldots&a_n\cr &\cr &\bigcirc\cr}\right).$$
 Since
 $$\eqalign{\max\left(\left\|\sum b^*_ib_i\right\|^{1/2}, \left\|\sum
 a_ia^*_i\right\|^{1/2}\right) &\le \sqrt n\, \max\{\|a_1\|,\ldots,
\|a_n\|,
 \|b_1\|,\ldots, \|b_n\|\}\cr
 &\le \sqrt n\, \max\left\{\left\|\sum a^*_ia_i\right\|^{1/2}, \left\|
 \sum
 b_ib^*_i\right\|^{1/2}\right\}}$$
 it follows that $\|R\t 1_A\|\le \sqrt n$. Hence $\|R\|_{cb} \le\sqrt
 n$ and
 thus
 $$\left\|{1\over 2}(I_{ S_n }+R)\right\|_{cb} \le {1\over 2}(1+\sqrt
 n).$$
 To prove the converse inequality it suffices to consider $n\ge 2$. Let
 $A$
 be the Cuntz algebra $O_n$ (cf. [C]), which is generated
by $n$ isometries
 $s_1,\ldots, s_n \in B(H)$ satisfying
 $$\leqalignno{s^*_is_j &= \delta_{ij}I&(4.1)\cr
 \sum^n_{i=1} s_is^*_i &= 1.&(4.2)}$$
 By (4.2) the element
 $$z = \left(\matrix{s^*_1\cr \vdots&\bigcirc\cr s^*_n\cr}\right)
 \oplus
 \left(\matrix{s_1&\ldots&s_n\cr &\cr &\bigcirc\cr}\right)$$
 in $ S_n \t A$ has norm $\|z\| = 1$, while
 $$\left({1\over 2}(I_{S_n}+R)\t I_A\right) (z) = {1\over
2}
 \left(\left(\matrix{s_1+s^*_1\cr \vdots &\bigcirc\cr
 s_n+s^*_n\cr}\right)
 \oplus \left(\matrix{s_1+s^*_1&\ldots&s_n+s^*_n\cr
 &\cr &\bigcirc\cr}\right)\right)$$
 has norm
 $$\eqalign{{1\over 2} \left\|\sum^n_{i=1} (s_i+s^*_i)^2\right\|^{1/2}
 &=
 {1\over 2} \sup \left\{
\sum^n_{i=1}\|(s_i+s^*_i)\xi\|^2 \mid \xi \in
 H,\ \|\xi\| = 1\right\}^{1/2}\cr
 &= {1\over 2} \sup \left\{\left(\left(\sum^n_{i=1} s^*_is_i + s_is^*_i
 +
 s^2_i+(s^*_i)^2\right)\xi,\xi\right) \ \bigg|\ \xi \in
 H,\ \|\xi\|=1\right\}^{1/2}.}$$
 By (4.1) and (4.2), $\sum\limits^n_{i=1} s^*_is_i = nI$ and
 $\sum\limits^n_{i=1} s_is^*_i = I$. Set
 $$v = {1\over \sqrt n} \sum^n_{i=1} s^2_i.$$
 By (4.1), $v^*v = I$ so $v$ is an isometry. By (4.1) the range of $v$
 is
 orthogonal to the range of the isometry $s_1s_2$:\ Indeed for
 $\xi,\eta \in
 H$
 $$\eqalign{(v\xi, s_1s_2\eta) &= {1\over \sqrt n}
\left(\sum_i
 (s^*_2 s^*_1 s^2_i  \xi,\eta)\right)\cr
 &= {1\over \sqrt n} (s^*_2 s_1\xi,\eta)\cr
 &= 0.}$$
 Hence $v$ is a non-unitary isometry. Therefore the point spectrum of
 $v^*$
 contains the open unit disk $D$ (cf. e.g. [KP], p.253).
Hence also
 the ``numerical range'' of $v$
 $$\{(v\xi,\xi)\mid \|\xi\| = 1\} = \{(\xi, v^*\xi)\mid \|\xi\| = 1\}$$
 contains the open unit disk. In particular the number 1 is in the
 closure
 of this set. Therefore
 $$\eqalign{&\sup_{\|\xi\|=1} \left(\left(\sum^n_{i=1} s^*_is_i +
 s_is^*_i +
 s^2_i+(s^*_i)^2\right)\xi,\xi\right)\cr
 &\quad =  n +1+2\sqrt n\, \sup_{\|\xi\|=1} (\hbox{Re}(v\xi,\xi))\cr
 &\ge n +1+2\sqrt n\cr
 &= (1+\sqrt n)^2.}$$
 Hence $\big\|\big({1\over 2}(I_{S_n}+R)\t I_A\big)
(z)\big\| \ge {1\over
 2}(1+\sqrt n)\|z\|$, which proves (a).

 \n (b)\ Let $Q$ be a projection from $ S_n $ onto $E_n$. Set
 $\widehat Q = Q R = RQ R$. Then $\widehat Q$ is also
 a projection from $ S_n $ to $E_n$. Let $i_m$ denote the identity on
 $M_m$
 and $t_m$ the transposition of $M_m$. Then
 $$\widehat Q \t i_m = (R\t t_m) ( Q\t i_m)(R\t
t_m).$$
 Since $t_n\t t_m$ can be identified with transposition on $M_{nm}$,
 $\|t_n\t t_m\|=1$. Hence by the definition of $R$,
 $$\|R\t t_m\| \le 1.$$
 Therefore
 $$\|\widehat  Q\t i_m\| \le \| Q \t i_m\|,\quad m\in {\bf N}$$
 and so $\|\widehat Q\|_{cb} \le \| Q\|_{cb}$, and therefore also
 $$\left\|{1\over 2}( Q +\widehat  Q)\right\|_{cb} \le
 \| Q\|_{cb}.$$
 But
 $${1\over 2}( Q + \widehat Q) =  Q\left({1\over
 2}(I_{ S_n }+R)\right)$$
 and since $ Q$ is the identity on $E_n$, which is the range of
 ${1\over 2}(I_{ S_n }+R)$, we have ${1\over 2}( Q+\widehat
Q) =
 {1\over 2}(I_{ S_n }+R)$. Thus
 $$\| Q\|_{cb} \ge\left\|{1\over 2}(I_{ S_n }+R)\right\| =
{1\over 2}
 (1+\sqrt n).$$
 If $\psi\colon \ M_n\oplus M_n {\buildrel \hbox{onto}\over
 \longrightarrow}
 E_n$ is a projection of norm 1, then from the above
 $$\|\psi\|_{cb} \ge \|\psi_{| S_n }\|_{cb} \ge {1\over 2} (1+\sqrt
 n),$$
 proving (b).\qed

 \n {\bf Proof of Theorem 4.1.} Let $I_{E_n} = vw$ be a factorization
 of
 $I_{E_n}$ through an injective von~Neumann algebra $A$. By the
 injectivity
 of $A$, $w$ can be extended to a linear map $\widetilde w\colon \ M_n
 \oplus M_n \to A$ such that $\|\widetilde w\|_{cb} \le \|w\|_{cb}$
(cf. [Pa] Theorem 7.2). Clearly
 $ Q = v\widetilde w$ is a projection of $M_n\oplus M_n$ onto $E_n$.
 Hence by (b) in the preceding lemma
 $$\|v\|_{cb} \|w\|_{cb} \ge \|v\|_{cb} \|\widetilde w \|_{cb} \ge
 \| Q\|_{cb} \ge {1\over 2}(1+\sqrt n).$$
 This proves the announced result when $A$ is an injective von~Neumann
 algebra. If $A$ is a nuclear $C^*$-algebra, and $I_{E_n} = vw$ as
 above, we
 can extend $v$ to a $\sigma(A^{**}, A^*)$-continuous linear map
 $\tilde
 v\colon \ A^{**}\to E_n$ such that $\|\tilde v\|_{cb} = \|v\|_{cb}$.
 Since
 $A^{**}$ is an injective von~Neumann algebra (cf. e.g.
[CE]), we are now
 reduced to the preceding case.\qed

 \n {\bf Remark 4.3.} The constant ${1\over 2}(1+\sqrt n)$ is best
 possible
 in   Theorem 4.1 :\ Namely let $A = M_n\oplus M_n$, let
$w\colon \ E_n\to M_n
 \oplus M_n$ be the inclusion map and define a projection $v\colon
 \ M_n
 \oplus M_n\to E_n$ by
 $$v(x\oplus y) = {1\over 2}(I_{ S_n }+R) (xp\oplus py), \quad
x\oplus y \in
 M_n\oplus M_n$$
 where $p = \left(\matrix{1&0&\ldots&0\cr 0\cr \vdots&&\bigcirc\cr
 0\cr}\right)$. Then clearly $vw = I_{E_n}$ and
 $$\|v\|_{cb}=\|v\|_{cb} \|w\|_{cb} \le \left\|{1\over 2}
(I_{ S_n }+S)\right\|_{cb} =
 {1\over 2}(1+\sqrt n),$$
which indeed shows that Theorem 4.1 is sharp.

 In connection with Lemma 4.2 (b), note that there is
obviously a projection $P\colon \ M_n \oplus M_n\to
 E_n$ with (ordinary) norm $\|P\|\le 1$ (simply take
$P=v$ with $v$ as in Remark 4.3). However, we will show
below that
 the projection constant of $E_n\t M_n$ in $(M_n\oplus
M_n)\t M_n$ goes to
 infinity when $n\to \infty$. To see this it is clearer to place the
 discussion in a broader context.

 Let $S\subset B(H)$ be a closed subspace. We define $\lambda(S)$
 (resp.
 $\lambda_{cb}(S)$, $\lambda_n(S)$) to be the infimum of the constants
 $\lambda$ such that there is a projection $P\colon \ B(H)\to S$
 satisfying
 $\|P\|\le \lambda$ (resp. $\|P\|_{cb} \le \lambda$,resp.
 $\|I_{M_n}\otimes
 P     \|_{M_n (B(H))\to M_n(S)} \le \lambda$). Then by
the extension
 theorem of $c.b.$ maps (cf. [W,Pa]), these constants are
invariants of
 the ``operator space'' structure of $S$. By this we mean that if
 $S_1\subset B(K)$ is another operator space which is completely
 isometric
 to $S$ (resp. such that for some constant  $\lambda$ there is an
 isomorphism
 $u\colon \ S\to S_1$ {with} $ \|u\|_{cb}
\|u^{-1}\|_{cb} \le
 \lambda)$
 then $\lambda(S_1) = \lambda(S), \lambda_{cb}(S_1) = \lambda_{cb}(S)$,
 $\lambda_n(S_1) = \lambda_n(S)$ (resp. ${1\over \lambda} \lambda(S)
 \le
 \lambda(S_1) \le \lambda\lambda(S)$ and similarly for the other
 constants).

 By a simple averaging argument we can prove

 \proclaim Proposition 4.4. Let $S\subset B(H)$ be a closed subspace.
 Consider $M_n(S) = M_n\t S\subset B(\ell^n_2(H))$. Then\medskip
 \item{(i)} $\lambda_n(S) = \lambda(M_n(S))$.
 \item{(ii)} If $S$ is $\sigma(B(H), B(H)_*)$-closed in $B(H)$ then
 \medskip
 $$\lambda_{cb}(S) = \sup_{n \ge 1} \lambda_n(S).\leqno (4.3)$$
 For any infinite dimensional Hilbert space $K$ we have
 $$\lambda_{cb}(S) \leq \lambda(B(K)\t S).\leqno (4.4)$$
  Moreover let $B(K)\overline{\t} S$ denote
 the weak-$*$ closure of $B(K){\t} S$ in
$B(K\t H)$. Then,
$$\lambda_{cb}(S)= \lambda(B(K)\overline{\t} S).$$

 \n {\bf Proof.} (i)\ The inequality $\lambda(M_n\t S) \le
 \lambda_n(S)$ is
 obvious, so we turn to the converse. Assume that there is a projection
 $$P\colon \ M_n\t B(H)\to M_n\t S$$
 with $\|P\| \le \lambda$.

 Let ${\cal U}_n$ be the group of all $n\times n$ unitary matrices.
 Consider
 then the group $G = {\cal U}_n \times {\cal U}_n$ equipped with its
 normalized Haar measure $m$. We will use the representation
 $$\pi\colon \ G\to B(M_n,M_n)$$
 defined by
 $$\pi(u,v) x = uxv^*.$$
 We can define an operator $\widetilde P\colon\ M_n \t B(H) \to M_n \t
 B(H)$
 by the following formula
 $$\widetilde P = \int(\pi (u,v)\t I_{B(H)}) P(\pi(u,v)\t
 I_{B(H)})^{-1}
 dm(u,v).\leqno (4.5)$$
 Note that  $\pi(u,v)$ leaves $M_n\t S$ invariant so that the range of
 $\widetilde P$ is included in $M_n\t S$ and $\widetilde P$ restricted
 to
 $M_n\t S$ is the identity, hence $\widetilde P$ is a projection from
 $M_n
 \t B(H)$ onto $M_n\t S$. Moreover, by Jensen's inequality (notice that
 $\pi(u,v)\t I_{B(H)}$ is an isometry on $M_n \t B(H)$) we have
 $$\|\widetilde P\| \le \|P\| \le \lambda.$$
 Furthermore, using the translation invariance of $m$ in (4.5) we find
 $$\forall\ (u_0,v_0)\in G \qquad \widetilde P(\pi(u_0,v_0) \t
 I_{B(H)}) =
 (\pi(u_0,v_0) \t I_{B(H)})\widetilde P,\leqno(4.6)$$
 so that $\widetilde P$ commutes
with $\pi(u_0,v_0) \t I_{B(H)}$.
 By well
 known facts   this implies that $\widetilde P$ is of the
form
 $$\widetilde P = I_{M_n}\t  Q$$
 for some operator $ Q$ which has to be a projection onto
$S$.  Indeed, since $M_n$ is spanned by
${\cal U}_n$, the above formula (4.6) is equivalent to:
For all $a,b$ in $M_n$ and for all $x$ in $M_n \otimes
B(H)$,
$$\w  ((a \otimes  1)x(b \otimes  1)) = (a \otimes
1)\w  (x) (b \otimes  1). \leqno(4.7)$$
 Let $(e_{ij})_{i,j=1,...,n}$ denote the matrix units in
$M_n.$ Set  $x = e_{ij} \otimes
 y,$ where $y$ is in $B(H)$ and $i,j$ are in
$\{1,...,n\}.$  Applying (4.7) to
 $a = 1-e_{ii}$ and $b = 1-e_{jj}$, one gets \nobreak
     $(1-e_{ii}) \w  (e_{ij} \otimes  y)(1-e_{jj}) =0$
{\it i.e.}  $\w  (e_{ij} \otimes  y) = e_{ij} \otimes  z $
for some $z$ in $B(H)$ depending
 on $y,i$ and $j$. However applying (4.7) again, this time
with  $a=e_{ki}$ and
 $b=e_{jl}$ it follows that $z$ is independent of $i$ and
$j$. Hence
$\w  = I_{M_n} \otimes  Q$,
 for some operator $ Q$ (which has to be a projection onto
$S$).
 Finally, we conclude
 $$\|I_{M_n} \t  Q\| = \|\widetilde P\| \le \lambda$$
 hence $\lambda_n(S) \le \lambda(M_n\t S)$.
 This proves (i).

 We now check (ii). Consider an arbitrary closed subspace $S\subset
 B(H)$
 and let $\overline S$ be the $\sigma(B(H), B(H)_*)$-closure of $S$. We
 claim that there is an operator $ Q\colon \ B(H)\to \overline S$ such
 that $ Q_{|S} = I_S$ and  $\| Q\|_{cb} \le \sup\limits_n
 \lambda_n(S)$.

 \n Let $\varepsilon_n>0$ be such that $\varepsilon_n\to
0$. For each $n$ there
 is a projection $P_n\colon\ B(H)\to S$ such that
 $$\|I_{M_n} \t P_n\|_{M_n(B(H))\to M_n(S)} \le (1+\varepsilon _n)
 \lambda_n(S).\leqno (4.8)$$
 Let ${\cal U}$ be a non-trivial ultrafilter on $\bf N$. For any
 bounded
 sequence $(\alpha_n)$ of real numbers (or for any
relatively compact sequence in a topological space) we will
denote simply by
 $\lim\limits_{\cal U}
 \alpha_n$ the limit of $\alpha_n$ when $n\to \infty$ along ${\cal U}$.
 For
 any $x$ in $B(H)$ let
 $$ Q(x) = \lim_{\cal U} P_n(x)$$
 where the limit is in the $\sigma(B(H), B(H)_*)$-sense. Observe that
 $\| Q\| \le \lim\limits_{\cal U}\|P_n\| \le \sup\limits_n
 \lambda_n(S)$. More generally for any integer $m\ge 1$ we clearly have
 $$\forall\ y \in M_m \t B(H)\qquad (I_{M_m}\t  Q)(y) =
\lim_{\cal
 U}(I_{M_m}\t P_n)(y)$$
 hence $\|I_{M_m}\t  Q\| \le \lim\limits_{\cal U}
\|I_{M_m}\t P_n\|$
 but when $n\ge m$ we have obviously
 $$\|I_{M_m}\t P_n\| \le \|I_{M_n}\t P_n\|$$
 hence by (4.8) we obtain
 $$\|I_{M_m}\t  Q\| \le \lim_{\cal U} (1+\varepsilon_n)
\lambda_n(S)
 \le \sup_n \lambda_n(S),$$
 so that $\| Q\|_{cb} \le \sup_n\lambda_n(S)$. Clearly $
Q(B(H))
 \subset \overline S$ and $ Q_{|S} = I_S$. This proves our claim and in
 the case $\overline S = S$ we obtain (4.3).
 (Note that $\lambda_{cb}(S) \ge \sup\limits_n \lambda_n(S)$ is
 trivial.)
We now turn to (4.4). We may clearly assume $K=\ell_2$. Recall that
there
 is obviously a completely contractive projection $\pi_n\colon \
 B(\ell_2)\to M_n$ (here $M_n$ is considered as a subspace of
 $B(\ell_2)$ in
 the usual way) hence
 $$\lambda_n(S) = \lambda(M_n\t S) \le \|\pi_n\| \lambda(B(\ell_2)\t S)
 =
 \lambda(B(\ell_2)\t S)$$
 which implies by (4.3)
 $$\lambda_{cb}(S) \le \lambda(B(\ell_2)\t S).$$
This
 concludes the proof of (4.4).

\n To prove the last
assertion, note that $M_n(S)$ is clearly contractively
complemented in $B(\ell_2)\overline{\t}S$ hence we have
$$\lambda_{cb}(S)\le   \sup_{n \ge 1} \lambda_n(S)=
\sup_{n \ge 1} \lambda(M_n(S))\le
\lambda(B(\ell_2)\overline{\t}S).$$
To prove the converse inequality, note that
$B(\ell_2)\overline{\t}B(H)$ can be identified with the
space of  matrices $a=(a_{ij})_{i,j\in \nat}$ which
 are bounded on $\ell_2(H)$, and
$B(\ell_2)\overline{\t}S$ can be identified with the
subspace formed by all matrices with entries in $S$. Then
 if $P$ is a
bounded projection from $B(H)$ onto $S$, defining
$$\w (a)=(P(a_{ij}))_{i,j\in \nat}$$
we obtain a projection from $B(\ell_2)\overline{\t}B(H)$
to $B(\ell_2)\overline{\t}S$ with
$\|\w  \|\le \|P\|_{cb}$. To check this last
estimate, observe that the norm of an element
 $a=(a_{ij})_{i,j\in \nat}$
  in $B(\ell_2)\overline{\t}B(H)$ is the supremum over $n$
of the norms in $M_n(B(H))$ of the matrices
$(a_{ij})_{i,j\le n}$.
This yields the last
assertion. \qed

\proclaim Corollary 4.5. Let $H,K$ be Hilbert spaces. Consider a
 completely isometric embedding $E_n\to B(H)$. Then, if $\dim K =
 \infty$,
 for any projection $P$ from $B(K)\t B(H)$ to $B(K)\t E_n$ we have
 $$\|P\| \ge {1\over 2}(\sqrt n+1).$$
 \ A fortiori the same holds for any projection $P$
from $B(K\t H)$
 onto $ B(K) \t E_n$.

 \n {\bf Proof.} By the preceding statement,
 this follows from Theorem~4.1.

\proclaim Corollary 4.6. Let $M\subset B(H)$ be a von
Neumann subalgebra such that  $M$ is isomorphic (as a
von Neumann algebra) to $M_n(M)$ for some integer $n\ge
2$. Then if there is a bounded linear projection from
$B(H)$ onto $M$, there is also a completely bounded one.

\p Note that if $M$ is isomorphic  to $M_n(M)$, then
obviously it is isomorphic to $M_n(M_n(M))=M_{n^2}(M)$,
and similarly to $M_{n^3}(M)$, and so on. Hence
this follows clearly from the first two
parts of Proposition 4.4 and the observation preceding
Proposition 4.4. \qed

In particular we have using [V1]
\proclaim Corollary 4.7. Let $M\subset B(H)$ be a von
Neumann subalgebra. If $M$ is isomorphic to
the von Neumann algebra
$VN(F_n)$ (resp. $VN(F_\infty)$) associated to the free
group with $n>1$ generators (resp. countably many
generators)  then there is no bounded linear projection
from $B(H)$ onto $M$.

\p
First note that $VN(F_n)$ trivially embeds into
$VN(F_\infty)$ as a subalgebra  which  is the range of  a
completely contractive projection. Therefore by
Proposition 1.3 and Theorem 4.1  there is no {\it
completely } bounded projection from  $B(H)$ onto $M$
if $M$ is isomorphic to $VN(F_\infty)$. By [V1]
$M_n(VN(F_\infty))$ is isomorphic to $VN(F_\infty)$ for
all $n$. Hence Corollary 4.7 for  $VN(F_\infty)$ follows
from the preceding corollary. To obtain the  case
of finitely many generators, recall  the well
known fact that $F_\infty$ can be embedded in $F_n$ for
all $n\ge 2$.
(If $a,b$ are two of the generators of $F_n$,
then it is easy to check, that
	$b, aba^{-1},..., a^{n}ba^{-n},...$
    are free generators of a subgroup isomorphic
 to $F_\infty$.) Therefore if $M=VN(F_n)$ for $n>1$,
then $VN(F_\infty)$ is isomorphic to a von Neumann
subalgebra $N\subset M$, and since  $M$
is a finite von Neumann algebra,  $N$ is the
range of a conditional expectation, hence there
is a
  bounded projection from $M$ onto $N$. Since
there is no bounded projection from $B(H)$ onto  $N$ by
the first part of the proof, a fortiori there cannot exist
a bounded projection from $B(H)$ onto  $M$. \qed

For two operator spaces $E$ and $F$ of the same finite dimension $n$,
one
can define the complete version of the Banach-Mazur
distance between $E$ and $F$ by
$$d_{cb}(E,F) = \inf\{\|u\|_{cb}, \|u^{-1}\|_{cb}\},$$
where the infimum is taken over all invertible linear maps $u$ from $E$
to $F$. By Proposition~1.3 it follows that
$$d_{cb}(E_n, \hbox{span}\{\lambda(g_i)\mid i=1,\ldots, n\}) \le 2$$
for all $n\in {\bf N}$. The next proposition shows that the same
inequality
holds if the unitary operators $\lambda(g_1),\ldots, \lambda(g_n)$ are
replaced by a semicircular or circular system of operators
in the sense of Voiculescu [V1].

\proclaim Proposition 4.8. Let $n\in {\bf N}$ and let $x_1,\ldots, x_n$
be
a semicircular or circular system of operators on a Hilbert space, then
the
map $u\colon \ E_n\to \hbox{span}\{x_1,\ldots, x_n\}$ given by
$$u\colon \ \sum^n_{k=1}c_k\delta_k\longrightarrow \sum^n_{k=1}
c_kx_k,\qquad
c_1\in {\bf C}$$
satisfies $\|u\|_{cb} \|u^{-1}\|_{cb} \le 2$.

\p Assume first that $x_1,\ldots, x_n$ is a semicircular
system of selfadjoint operators in the sense of [V1]. By
[V2], we can exchange $x_1,\ldots, x_n$ with the operators
$$x_k = {1\over 2 } (s_k+s^*_k),\qquad k=1,\ldots, n$$
where $s_1,\ldots, s_n$ are the ``creation operators'' $\xi\to
e_i\otimes
\xi$ on the full Fock space
$${\cal H} = {\bf C} \otimes \left(\bigoplus ^\infty_{n=1} H^{\otimes
n}\right)$$
based on a Hilbert space $H$ with orthonormal basis $(e_1,\ldots,
e_n)$. In
particular $s_1,\ldots, s_n$ are $n$ isometries with orthogonal ranges,
and
therefore
$$\sum^n_{k=1} s_ks^*_k \le 1.$$
Hence, as in the proof of Proposition~1.1, we get that for any
$n$-tuple
$a_1,\ldots, a_n$ of elements in a $C^*$-algebra $A$,
$$\eqalign{\left\|\sum_k x_k\otimes a_k\right\| &\le {1\over 2}
\left(\left\| \sum_k s_k\otimes a_k\right\| + \left\|\sum_k
s^*_k\otimes
a_k\right\|\right)\cr
&\le {1\over 2} \left(\left\| \sum_k s_ks^*_k\right\|^{1/2}
\left\|\sum_k
a^*_ka_k\right\|^{1/2} + \left\|\sum_k s^*_ks_k\right\|^{1/2}
\left\|\sum_k
a_ka^*_k\right\|^{1/2}\right)\cr
&\le \max\left\{\left\|\sum_k a^*_ka_k\right\|^{1/2}, \left\|\sum_k
a_ka^*_k\right\|^{1/2}\right\}.}$$
Hence $\|u\|_{cb} \le 1$. To prove that $\|u^{-1}\|_{cb} \le 2$, notice
that by [V1], [V2], the $C^*$-algebra generated by
$x_1,\ldots, x_n$ and 1 has a trace
$$\tau \colon \ C^*(x_1,\ldots, x_n,1) \to {\bf C}$$
(namely the vector-state given by a unit vector in the ${\bf C}$-part
of
the Fock space ${\cal H}$), with the properties:
$$\tau(1) = 1, \quad \tau(x^2_k) = {1\over 4}\quad
\hbox{and}\quad \tau(x_kx_\ell) = 0\qquad k\ne \ell.$$
Let $a_1,\ldots, a_n$ be $n$ operators in a $C^*$-algebra $A$, and let
$S(A)$ denote the state space of $A$. Then
$$\eqalign{\left\|\sum_k x_k\otimes a_k\right\|^2 &\ge \sup_{\omega\in
S(A)} (\tau\otimes\omega) \left(\left( \sum_k x_k \otimes a_k\right)^*
\left(\sum_\ell x_\ell \otimes a_\ell\right)\right)\cr
&= {1\over 4} \sup_{\omega\in S(A)} \omega\left(\sum_k
a^*_ka_k\right)\cr
&= {1\over 4} \left\|\sum_k a^*_ka_k\right\|}$$
and similarly $\left\|\sum\limits_k x_k\otimes a_k\right\|^2 \ge
{1\over 4}
\left\|\sum\limits_k a_ka^*_k\right\|$. Hence
$$\left\|\sum_k x_k\otimes a_k\right\| \ge {1\over 2}
\max\left\{\left\|
\sum_k a^*_ka_k\right\|^{1/2}, \left\|\sum_k
a_ka^*_k\right\|^{1/2}\right\}$$
proving $\|u^{-1}\|_{cb} \le 2$.

Assume finally that $y_1,\ldots, y_n$ is a circular system. Then
$$y_k = {1\over \sqrt 2} (x_{2k-1} + ix_{2k}),\qquad k=1,\ldots, n,$$
where $(x_1,\ldots, x_{2n})$ is a semicircular system of selfadjoint
operators. Therefore the statement about circular systems in
Proposition~4.8 follows from the one on semicircular systems by
observing,
that the map
$$\sum^n_{k=1} c_k\delta_k \to {1\over \sqrt 2} \sum^n_{k=1}
c_k(e_{2n-1} +
e_{2k})$$
defines a $cb$-isometry of $E_n$ onto its range in $E_{2n}$.
\qed

To conclude this paper we give a generalization of Proposition 1.1 to
free products of discrete
groups, or more generally free products of $C^*$-probability spaces in
the
sense of [V1] and [V2]. We refer to [V1] and [V2] for the terminology.

\proclaim { Proposition 4.9}. Let $(A,\ph)$ be a $C^*$-algebra equipped
with
a faithful state $\ph$. Let $(A_i)_{i \in I}$ be a free family of
unital
$C^*$-subalgebras of $A$ in the sense of [V1] or [V2]. Consider
elements $x_i
\in A_i$ such that for some $\delta > 0$
$$\forall~-i \in I~~~~~\N {x_i} \leq 1,\ \ph(x_i)=0\  {\rm
and}~-\min \{\ph(x_i^* x_i),\ph(x_i x_i^*)\} \geq
\delta^2$$ Then, for all finitely supported families
$(a_i)_{i \in I}$ in $B(H)$ ($H$ Hilbert) we have
$$\delta \max \{\N {\sum a_i^* a_i}^{1/2}, \N {\sum a_i a_i^*}^{1/2}\}
\leq \N {\sum x_i \otimes a_i}
\leq 2 \max \{\N {\sum a_i^* a_i}^{1/2},\N {\sum a_i a_i^*}^{1/2}\}.
\leqno
(4.9)$$

\n {\bf Proof.}\  We may assume that $I$ is finite. The lower bound in
(4.9) is
proved exactly as in the semicircular case. To prove the upper bound we
will
prove that $A$ can be faithfully represented as a $C^*$-algebra of
operators
on a Hilbert space $H$, such that $x_i$ admits a decomposition $x_i =
u_i +
v_i$ with $u_i,v_i$ in $B(H)$ and
$$\N {\sum u_i^* u_i} \leq 1~-{\rm and}~-\N{\sum v_i v_i^*} \leq
1.\leqno
(4.10)$$
The upper bound in (4.9) then follows as in the semicircular case.

\n Following the notation of [V 2, pp. 558-559], we let $(H_i,\xi_i)$
be the
space of the GNS-representation $\pi_i = \pi_{\ph|A_i}$. In particular
$\xi_i$
is a unit-vector in $H_i$ and
$$\ph(x) = (\pi_i(x) \xi_i,\xi_i) \ {\rm when}\   x \in
A_i.$$ Then $A$ can be realized as the $C^*$-algebra of
operators on the Hilbert space $$(H,\xi) = *_{i \in
I}-(H_i,\xi_i)$$ generated by  $\bigcup \limits_{i \in I}
\lambda_i\rond\pi_i(A_i)$, where $\lambda_i : B(H_i) \to
B(H)$ is the $*$-representation defined in [V2, sect.
1.2]. For simplicity of notation we will identify $A_i$
with its range in $B(H)$, \ie we set  $$\lambda_i \rond
\pi_i(x) = x\ {\rm when}\ x \in A_i.$$ Let $x \in A_i$.
Corresponding to the decomposition $$H_i = H_i^0 \oplus
\comp {\xi_i}.$$ we can write $\pi_i(x)$ as a $2 \times
2$ matrix $$\pi_i(x) = \left ( \matrix{
b & \eta\cr
\zeta^* & t\cr
}\right)$$
where $b \in B(H_i^0), \eta,\zeta \in H_i^0$ and $t \in \comp$.
(Here we identify  $\eta,\zeta$ with the corresponding
linear maps from $\comp$ to $H_i^0$, and
  we also   identify  $\comp$ with
$\comp{\xi_i}$.)
 The action of $x
= \lambda_i \rond \pi_i(x)$ on $*_{i \in I} (H_i,\xi_i)$ can now be
explicitly
computed from [V 2, sect. 1.2]. One finds :

 (4.11)~~~~~ $x \xi = \eta \otimes\xi  + t \xi$,

(4.12)~~~~~ $x(h_1\otimes...\otimes h_n) = bh_1\otimes...\otimes h_n +
(h_1,\zeta) h_2 \otimes...\otimes h_n$ when $n \geq 1$,\hfill$\ h_k \in
H_{i_k}^0,\  i = i_1 \not= i_2\not=...\not=i_n$,

(4.13)~~~~~ $x(h_1\otimes...\otimes h_n) = \eta \otimes h_1
\otimes...\otimes h_n
+ th_1 \otimes...\otimes h_n$, when $n \geq 1$,$\ h_k \in
H_{i_k}^0,$\hfill $\ i\not=i_1\not=i_2\not=...\not=i_n$

\noindent where $h_2\otimes...\otimes h_n = \xi$ for $n=1$.

\n  Let $e_i \in B(H)$ be
the orthogonal projection of $H$ onto the subspace
$$H_i = \oplus_{n=1}^\infty (\oplus (H_{i_1} \otimes...\otimes
H_{i_n}))$$
where the second direct sum contains all $n$-tuples $(i_1,...,i_n)$ for
which
$i = i_1\not=i_2\not=...\not=i_n$. From (4.11), (4.12) and (4.13) one
gets
for all $  x$ in $A_{i}$
$$(1 -e_i) x(1 - e_i) = \ph(x) (1 - e_i) \leqno (4.14)$$
where we have used that
$$t = (\pi_i(x) \xi_i,\xi_i) = \ph(x).$$
Let now $x_i \in A_i,\  \N {x_i}-\leq 1,\  \ph(x_i) = 0$. Then by
(4.14)
$$(1 - e_i) x_i(1 - e_i) = 0$$
Thus $x_i = u_i + v_i$, where
$$u_i = x_i e_i~-{\rm and}~-v_i = e_i x_i (1 - e_i).$$
Since $\N {x_i}-\leq 1$, and since $(e_i)_{i \in I}$ is a set of
pairwise
orthogonal projections,
$$\sum_{i \in I} u_i^* u_i \leq \sum_{i \in I} e_i \leq 1$$
and
$$\sum_{i \in I} v_i v_i^* \leq \sum_{i \in I} e_i \leq 1.$$
This completes the proof of proposition 4.9. \qed
\vfill\eject

\centerline {\bf REFERENCES}
\def\bg{\bigskip\goodbreak}
\pagetitrefalse
\bg
\bg
\bg \bg \bg

\ref BP & Blecher D. and  Paulsen V.& Tensor products
of operator spaces.&
 J. Funct. Anal. &99 (1991) 262-292.

\ref C &  Cuntz J. & Simple C*-algebras generated by
isometries.& Comm. Math. Phys.&
       57 (1977),173-185.

\ref CE & Choi M.D. and Effros E.G. &Nuclear C*-algebras
and injectivity. The
       general case.& Indiana Univ. Math. J.& 26 (1977),
443-446.

\ref CS & Christensen E. and Sinclair A. & On von
Neumann algebras which are complemented subspaces of
$B(H)$.& Preprint.&

\ref DCH & de Canni\`ere J. and Haagerup U. & Multipliers of the
Fourier algebras
of some simple Lie groups and their discrete subgroups.& Amer. J. Math.
& 107
(1985), 455-500.

 \ref {ER} &Effros E. and Ruan Z.J.& A new approach to
operators spaces.&
 Canadian Math. Bull.&
34 (1991) 329-337.

\ref E & Eymard P.& L'alg\`ebre de Fourier d'un groupe
localement compact. & Bull. Soc. Math.
France   & 92 (1964), 181-236.

\ref G & Grothendieck A.& R\'esum\'e de la th\'eorie m\'etrique des
produits tensoriels
topologiques. & Boll.. Soc. Mat. S$\tilde{a}$o-Paulo & 8 (1956), 1-79.

\ref H1 & Haagerup U. & Solution of the similarity
problem for cyclic representations of $C^*$-algebras. &
Annals of Math. & 118 (1983), 215-240.

\ref H2 & Haagerup U. & An example of a non-nuclear
$C^*$-algebra which has the metric approximation property.
& Inventiones Mat. & 50 (1979), 279-293.

\ref KP & Kadison R.V. and Petersen G.K. &Means and
convex combinations of
 of unitary operators,& Math. Scand.  &57 (1985),
249-266.

\ref KR & Kadison R. and Ringrose J. & Fundamentals of
the theory of operator algebras, Vol. II  & Academic
Press, New-York    & 1986.

\ref K & Kwapie\'n S. & Decoupling inequalities for
polynomial chaos. & Ann. Probab.  &15 (1987) 1062-1071.

\ref LPP & Lust-Piquard F. and Pisier G. & Non commutative Khintchine
and Paley
inequalities. & Arkiv fUr Mat. & 29 (1991), 241-260.

\ref Pa & Paulsen V. & Completely bounded maps on
$C^*$-algebras and invariant operator ranges.& Proc.
Amer. Math. Soc. & 86 (1982) 91-96.

\ref P1 & Pisier G. && Grothendieck's theorem for
noncommutative $C^*$-algebras with an appendix on
Grothendieck's constants. &J. Funct. Anal. 29 (1978)
397-415.

\ref P2 & Pisier G. & Random series of trace class
operators. & Proceedings IV\up{o} C.L.A. P.E.M. (Mexico
Sept. 90). & To appear.

\ref P3 & Pisier G. && Factorization of linear operators
and the Geometry of Banach spaces.& CBMS (Regional
conferences of the A.M.S.) n\up o 60, (1986), Reprinted
with corrections 1987.

\ref P4 & Pisier G. &  Espace de Hilbert d'op\'erateurs
et interpolation complexe. & C. R. Acad. Sci. Paris &
 S\'erie I, 316 (1993) 47-52.

\ref V1 & Voiculescu D. & Circular and semicircular
systems and free product factors. &In Operator algebras,
unitary representations, Enveloping algebras, and
invariant theory, (edited by A.Connes, M.Duflo,
A.Joseph, et R. Rentschler) Colloque en l'honneur de
J.Dixmier. Birkhauser, Progress in Mathematics &
vol. 92 (1990) 45-60.

\ref V2 & Voiculescu D. & Symmetries of some reduced
free product $C^*$-algebras, in & "Operator algebras and
their connections with Topology and Ergodic Theory",
Springer Lecture Notes in Math. &1132
(1985) 556-588.

\ref W & Wittstock G. & Ein operatorwertiger Hahn-Banach
satz,& J. Funct. Anal.& 40 (1981) 127-150.

\vskip12pt

Addresses

U.Haagerup:

Odense University

DK 5230 Odense, DENMARK
\vskip12pt
G.Pisier:

Texas A. and M. University

College Station, TX 77843, U. S. A.

and

Universit\'e Paris 6

Equipe d'Analyse, Bo\^\i te 186,

75255 Paris Cedex 05, France

\end